\documentclass[preprint,12pt]{elsarticle}

\usepackage{amssymb}
\usepackage{amsmath}
\usepackage{amsthm,amssymb}
\usepackage{parskip}
\usepackage{multirow}
\usepackage[utf8]{inputenc}
\usepackage[english]{babel}
\newtheorem{theorem}{Theorem}[section]

\newtheorem{lemma}[theorem]{Lemma}

\usepackage{pdfpages}
\usepackage{graphicx}
\usepackage{graphics}
\usepackage {framed}
\usepackage {amsxtra}
\usepackage {enumerate}
\usepackage{commath}
\usepackage[margin=1 in]{geometry}
\usepackage{float}
\usepackage[caption = false]{subfig}
\journal{arXiv}

\begin{document}

\begin{frontmatter}



\title{Weak-consistent dynamic correlation estimators for\\ Brownian motion pairs and for Geometric Brownian motion pairs}



\author{Majnu John$^{a,b,c}$\footnote{Corresponding author: Department of Mathematics, $308$ Roosevelt Hall, $130$ Hofstra University, Hempstead, NY 11549.  e-mail: {\sf Majnu.John@hofstra.edu}}, Yihren Wu$^{d}$}

\address{$^{a}$Departments of Mathematics and of Psychiatry, \\Hofstra University, Hempstead, NY.}
\address{$^{b}$The Feinstein Institutes of Medical Research, \\NorthWell Health System Manhasset, NY.}
\address{$^{c}$Division of Psychiatry Research, \\Zucker Hillside Hospital, Glen Oaks, NY.}
\address{$^{d}$Department of Mathematics, \\Hofstra University, Hempstead, NY.}

\begin{abstract}

Estimating dynamic correlation between a pair of time series is of importance in many applications. We present new estimators for the dynamic correlation between a pair of correlated Brownian motions and separately for dynamic correlation between a pair of correlated Geometric Brownian motions. We show that, as the sample size increases, all estimators presented in this paper converge in probability to the underlying true dynamic correlation.

\end{abstract}

\begin{keyword}

dynamic correlation, correlated bivariate Brownian motion, correlated bivariate Geometric Brownian motion, weak consistency




\end{keyword}

\end{frontmatter}





\noindent

\section{Introduction}

Dynamic correlation between a pair of time series is of long-standing interest in many fields such as economics [1,2] and neuroscience [3]. Recently, a nonparametric estimator of dynamic correlation, robust in the presence of extreme values, was introduced in [4]. The new estimator may be conceptualized based on a particular algorithm for conversion of time series into a weighted graph and considering the edge-weights of the graph. In this paper we consider extensions of the estimators presented in [4]. Specifically, we present estimators for a pair of correlated bivariate Brownian motions and estimators for a pair of correlated Geometric Brownian motions. The main focus of the paper is in showing weak consistency of the estimators.

The paper is structured as follows. The new estimators are presented in section 2. Next, in section 3, the weak consistency results are presented. Conclusions are summarized in the last section.

\section{New estimators}

First we focus on a correlated (bivariate) Brownian motion pair denoted by ($X_{t}, Y_{t}$), $t \in$ [$0, T$], and present a weakly-consistent estimator for the dynamic correlation, $\rho_{t} = \gamma_{t}/\sigma_{t}^{x}\sigma_{t}^{y} = \gamma_{t}/t$ between $X_{t}$ and $Y_{t}$. Assume $n$ points are sampled from each time series at time points $\{t_{1}, \ldots, t_{n} = T \} \equiv I$. We assume without loss of generality, $\mathbb{E}(X_{t}) = 0 = \mathbb{E}(Y_{t}), \forall t$ and for convenience, we take $t_{i} = i, \forall i$; however, our results will hold in general, even without these simplifications.
The estimators that we consider are $\hat{\rho}_{u}^{\,q,p} = \hat{\gamma}_{u}^{q,p}/(\hat{\sigma}_{u}^{x, q,p}\hat{\sigma}_{u}^{y, q,p})$, $u \in I, q, p \geq 0$, where \begin{equation}\label{eq1} \hat{\gamma}_{u}^{\,q,p} = \frac{1}{T-1}\sum\limits_{\substack{v \in I \\ v \neq u}} \frac{(v^{q}X_{u} - v^{-p}X_{v})(v^{q}Y_{u} - v^{-p}Y_{v})}{(u-v)^{2}},\end{equation} \[ \left(\hat{\sigma}_{u}^{x,\,q,p}\right)^{2} = \frac{1}{T-1}\sum\limits_{\substack{v \in I \\ v \neq u}} \frac{(v^{q}X_{u} - v^{-p}X_{v})^{2}}{(u-v)^{2}},\;\;\; \left(\hat{\sigma}_{u}^{y,\,q,p}\right)^{2} = \frac{1}{T-1}\sum\limits_{\substack{v \in I \\ v \neq u}} \frac{(v^{q}Y_{u} - v^{-p}Y_{v})^{2}}{(u-v)^{2}}.\]

The estimators that we originally considered in [4] were special cases of the estimator presented in Eq. (\ref{eq1}), specifically with $p = q = 0$. The main reason that we consider the more general version is that, while studying the weak consistency of the original estimators, we realized that the asymptotic bias is non-zero for those estimators. In other words, weak consistency does not hold in the special case $p = q = 0$. In the general case, we give in section 3 below, the specific range of values for $p$ and $q$ for which weak consistency can be shown to hold.

Next we consider a correlated (bivariate) Geometric Brownian motion pair ($R_{t}, S_{t}$), $t \in [0,T]$. The most general form that we should consider is \[ R_{t} = R_{0}\exp\{(\mu_{R} - \frac{\sigma_{R}^{2}}{2})t + \sigma_{R}W_{t} \}\; \mathrm{and}\;S_{t} = S_{0}\exp\{(\mu_{S} - \frac{\sigma_{S}^{2}}{2})t + \sigma_{S}U_{t} \},  \] where ($W_{t}, U_{t}$) is a correlated Brownian motion pair. For ease of exposition, we restrict our attention to the case with $R_{0} = S_{0} = 1$ and $\mu_{R} = \mu_{S} = \sigma_{R}^{2}/2 = \sigma_{S}^{2}/2$ and $\sigma_{R} = \sigma_{S} = \sigma$ so that \[ R_{t} = e^{\sigma W_{t}}\; \mathrm{and}\; S_{t} = e^{\sigma U_{t}} \] is the Geometric Brownian motion pair that we consider. Note that in this case we have $\mathbb{E}(R_{t}) = \mathbb{E}(S_{t}) = e^{\sigma^{2}t/2}$ and $\mathrm{Var}(R_{t}) = \mathrm{Var}(S_{t}) = e^{\sigma^{2}t}(e^{\sigma^{2}t} - 1)$. We also mention here that we use the same notation ($\rho_{t}$ and $\gamma_{t}$, respectively) for the correlation and covariance between $R_{t}$ and $S_{t}$ as well. That is, whether $\rho_{t}$ (similarly $\gamma_{t}$) stands for the correlation (covariance) between the Brownian motion pair ($X_{t}, Y_{t}$) or for the Geometric Brownian motion pair ($R_{t}, S_{t}$) will hopefully be understood from the context. We consider two estimators for the dynamic correlation $\rho_{t} = \gamma_{t}/\sigma_{t}^{R}\sigma_{t}^{S} = \gamma_{t}/(\sigma_{t,W}\sigma_{t,U})$ between $R_{t}$ and $S_{t}$. In the first case \begin{align} \displaystyle \hat{\gamma}_{t} = \frac{1}{e^{c\sigma^{2}T}}\sum_{k = 1}^{T}& \left\{ \left[e^{\frac{-b\sigma^{2}k}{2}}\left(e^{\sigma W_{k}}-e^{\frac{\sigma^{2}k}{2}}\right) - e^{\frac{a\sigma^{2}k}{2}}\left(e^{\sigma W_{t}}-e^{\frac{\sigma^{2}t}{2}}\right) \right] \right. \nonumber \\
& \left. \times \left[e^{\frac{-b\sigma^{2}k}{2}}\left(e^{\sigma U_{k}}-e^{\frac{\sigma^{2}k}{2}}\right) - e^{\frac{a\sigma^{2}k}{2}}\left(e^{\sigma U_{t}}-e^{\frac{\sigma^{2}t}{2}}\right) \right] \right\}.
\end{align}

In the second case, we have the estimator for $\gamma_{t}$ as
\begin{equation}\displaystyle \hat{\gamma}_{t} = e^{-c\sigma^{2}T}\sum_{k=1}^{T} \left\{e^{a\sigma^{2}k}[(e^{\sigma W_{t}} - e^{\sigma^{2}t/2})(e^{\sigma U_{t}} - e^{\sigma^{2}t/2}) ] - e^{-b\sigma^{2}k}[(e^{\sigma W_{k}} - e^{\sigma^{2}k/2})(e^{\sigma U_{k}} - e^{\sigma^{2}k/2})  ] \right\} \end{equation} By considering $W_{k} = U_{k}$ for all $k$ (and $W_{t} = U_{t}$) in equations (2) and (3) we get the corresponding estimates for $\sigma_{t,W}^{2}$ and $\sigma_{t,U}^{2}$. Here $a, b$ and $c$ are real-valued constants with appropriate ranges given in section 3.

If the correlation between $W_{t}$ and $U_{t}$ is denoted by $r_{t}$, then it is easy to verify that \[r_{t} = \frac{1}{\sigma^{2}t}\log \left[1 + \rho_{t}(e^{\sigma^{2}t}-1) \right]\;\mathrm{or}\; \rho_{t} = \frac{e^{r_{t}\sigma^{2}t} - 1}{e^{\sigma^{2}t} - 1 }. \]Based on this, we may also consider an estimate for $\rho_{t}$ by plugging in an estimate for $r_{t}$ (for example, the one given in Eq. (1)). However we do not pursue such estimators in this paper as our focus is on estimators that are generalizations of the original estimators that we considered in [4].

\section{Results}

\subsection{Weak consistency of dynamic correlation estimators in the Brownian Motion case}

Our main result in this subsection is the following theorem. \\

\begin{theorem} $\hat{\rho}_{u}^{\,q,p} \rightarrow \rho_{u}$ in probability for each $u$ and $p > q = 1/2$, as $T \rightarrow \infty$. \end{theorem} We first make a few remarks which will be used later.

\textbf{Remark 1:} A random variable $V$ follows a variance-gamma distribution with parameters $r >0$, $\theta \in \mathbb{R}$, $\sigma > 0$, $\mu \in \mathbb{R}$ (denoted VG($r, \theta, \sigma, \mu$)) if it has probability density function given by
\[ p_{VG}(x; r, \theta, \sigma, \mu) = \frac{1}{\sigma\sqrt{\pi}\Gamma(\frac{r}{2})} \exp \left(\frac{\theta}{\sigma^{2}}(x-\mu) \right)\left(\frac{\abs{x-\mu}}{2\sqrt{\theta^{2} + \sigma^{2}}} \right)^{\frac{r-1}{2}}K_{\frac{r-1}{2}}\left(\frac{\sqrt{\theta^{2} + \sigma^{2}}}{\sigma^{2}}\abs{x-\mu} \right), \] where $x \in \mathbb{R}$ and $K_{\nu}(\cdot)$ is a modified Bessel function of the second kind [5,6]. We have $\mathbb{E}(V) = \mu + r\theta$, $\mathrm{Var}(V) = r(\sigma^{2} + 2\theta^{2})$ [7].

\textbf{Remark 2:} If ($X,Y$) denote a bivariate normal random vector with zero means, variances ($\sigma_{x}^{2}, \sigma_{x}^{2}$) and correlation coefficient $\rho$, then $Z = XY$ follows VG($1, \rho\sigma_{x}\sigma_{y}, \sigma_{x}\sigma_{y}\sqrt{1-\rho^{2}}, 0$) [5]. In particular, based on Remark 1, $\mathbb{E}(Z) = \rho\sigma_{x}\sigma_{y},\; \mathrm{Var}(Z) = (1 + \rho^{2})\sigma_{x}^{2}\sigma_{y}^{2}\; \mathrm{and}\; \mathbb{E}(Z^{2}) = (1 + 2\rho^{2})\sigma_{x}^{2}\sigma_{y}^{2}$, which could also be obtained using direct calculation with bivariate density.

\textbf{Remark 3:} Based on Remark 2, for the bivariate Brownian motion ($X_{t}, Y_{t}$) mentioned in the first paragraph, we have $\mathbb{E}(X_{t}Y_{t}) = t\rho_{t}$, $\mathrm{Var}(X_{t}Y_{t}/t^{2}) = (1+\rho_{t}^{2})/t^{2} \leq 2/t^{2}$, $\mathbb{E}(X_{t}Y_{t})^{2}/t^{4} = [(1+\rho_{t}^{2})t^{2} + \rho_{t}^{2}t^{2}]/t^{4} \leq 3t^{-2}$ and for $s < t$, \begin{align} \mathbb{E}\left(\frac{X_{t}Y_{t}X_{s}Y_{s}}{t^{2}s^{2}}\right) &= \frac{1}{t^{2}s^{2}}\{\mathbb{E}[(X_{t}-X_{s})(Y_{t}-Y_{s})X_{s}Y_{s}] + \mathbb{E}[(Y_{t}-Y_{s})X_{s}^{2}Y_{s}] \nonumber \\ &\;\;\;\;\;\;\;\;\;\;\;\;\;\;+ \mathbb{E}[(X_{t}-X_{s})X_{s}Y_{s}^{2}] + \mathbb{E}[X_{s}^{2}Y_{s}^{2}]\} \nonumber \\ & \leq \frac{C_{1}(t-s)s + C_{2}(t-s)^{1/2}s^{3/2} +  C_{3}s^{2}}{t^{2}s^{2}} \label{ineqRmk3pf} \\ &\leq C \left(\frac{1}{ts} + \frac{\sqrt{1 - (s/t)}}{t^{3/2}s^{1/2}} + \frac{1}{t^{2}} \right) \leq \frac{C}{s^{2}}, \;\mathrm{since}\; 1 \leq s < t. \nonumber \end{align} Here $C_{1}, C_{2}, C_{3}$ and $C$ are generic finite positive constants.

In order to get the inequality in Eq. (\ref{ineqRmk3pf}), we first note that $X_{t}$ and $Y_{t}$ may be written as \begin{equation} \label{eqXY1} X_{t} = \rho_{t}Y_{t} + \sqrt{1 - \rho_{t}^{2}}X_{t}^{(1)}, \; Y_{t} = \rho_{t}X_{t} + \sqrt{1 - \rho_{t}^{2}}Y_{t}^{(1)},  \end{equation} where $\{X_{t}^{(1)}\}$ and $\{Y_{t}^{(1)}\}$ are mean zero Brownian motions with $\{X_{t}^{(1)}\}$ independent of $\{Y_{t}\}$ and $\{Y_{t}^{(1)}\}$ independent of $\{X_{t}\}$. Using Eq. ($\ref{eqXY1}$) we may write \begin{align} \label{XXminusX} (X_{t}-X_{s})X_{s} = (\rho_{t}\rho_{s})(Y_{t} - Y_{s})Y_{s} + (\rho_{t} - \rho_{s})Y_{s}^{2} &+ \sqrt{1-\rho_{t}^{2}}(X_{t}^{(1)} - X_{s}^{(1)})Y_{s} \\ & + [\sqrt{1-\rho_{t}^{2}} - \sqrt{1-\rho_{s}^{2}} ]X_{s}^{(1)}Y_{s}. \nonumber\end{align} so that \begin{align} \label{Rmk3term1ineq} (X_{t}-X_{s})(Y_{t}-Y_{s})X_{s}Y_{s}  \leq (Y_{t} - Y_{s})^{2}Y_{s}^{2} + 2|Y_{t} - Y_{s}||Y_{s}|^{3} &+ 2|X_{t}^{(1)} - X_{s}^{(1)}||Y_{t} - Y_{s}|Y_{s}^{2} \nonumber \\ &+2|X_{s}^{(1)}||Y_{t}-Y_{s}|Y_{s}^{2}. \nonumber \end{align} Hence, using the fact that for a Brownian motion $\{B_{u}\}$, $\mathbb{E}|B_{u}|^{\alpha} \leq K u^{\alpha/2}$, \[\mathbb{E}[(X_{t}-X_{s})(Y_{t}-Y_{s})X_{s}Y_{s}] \leq C_{1}(t-s)s + C_{2}\sqrt{(t-s)s^{3}}. \] Using Eq. (\ref{XXminusX}) again, we will get \[\mathbb{E}[(X_{t}-X_{s})X_{s}Y_{s}^{2}] \leq C_{2}\sqrt{(t-s)s^{3}}\;\mathrm{and}\;\mathrm{similarly}\; \mathbb{E}[(Y_{t}-Y_{s})X_{s}^{2}Y_{s}] \leq C_{2}\sqrt{(t-s)s^{3}}.\] We also have $\mathbb{E}(X_{s}^{2}Y_{s}^{2}) \leq 4s^{2}$. Putting all these together we get the inequality in Eq. (\ref{ineqRmk3pf}). This completes Remark 3.

We denote $S_{a} = \sum_{n=1}^{\infty} n^{-a}$ for $a > 0$; $S_{a} < \infty$ for $a > 1$.  Although the primary focus of this subsection and the main result is for the case $p > q = 1/2$, we have the following result for $p = q = 0$.\\

\begin{lemma} $\mathrm{Var}(\hat{\gamma}_{i}^{0,0})$ $\rightarrow 0$, $\mathrm{Var}([\hat{\sigma}_{i}^{x,0,0}]^{2})$ $\rightarrow 0$ and $\mathrm{Var}([\hat{\sigma}_{i}^{y,0,0}]^{2})$ $\rightarrow 0$, for each $i$, as $T \rightarrow \infty$. \end{lemma}

\textbf{Proof:}  \begin{align} \mathrm{Var}(\hat{\gamma}_{i}^{0,0}) \leq \mathbb{E}(\hat{\gamma}_{i}^{0,0})^{2} &= \frac{1}{(T-1)^{2}}\sum\limits_{\substack{j=1 \\ j \neq i}}^{T}\frac{\mathbb{E}[(X_{i}-X_{j})(Y_{i}-Y_{j})]^{2}}{(i-j)^{4}} \nonumber \\ &+ \frac{1}{(T-1)^{2}}\sum\limits_{\substack{j=1 \\ j \neq i}}^{T}\sum\limits_{\substack{k=1 \\ k \neq i,j}}^{T}\frac{\mathbb{E}[(X_{i}-X_{j})(Y_{i}-Y_{j})(X_{i}-X_{k})(Y_{i}-Y_{k})]}{(i-j)^{2}(i-k)^{2}}  \label{lmpq0pf1} \end{align} Using one of the facts mentioned in Remark 3 in the inequality below, \[ 1^{st}\; \mathrm{term}\; = \frac{1}{(T-1)^{2}} \left\{\sum_{j=1}^{i-1}\frac{\mathbb{E}(X_{i-j}Y_{i-j})^{2} }{(i-j)^{4}} + \sum_{j=i+1}^{T}\frac{\mathbb{E}(X_{j-i}Y_{j-i})^{2} }{(j-i)^{4}} \right\} \leq \frac{CS_{2}}{(T-1)^{2}} = O(T^{-2}), \] where $C$ is a finite positive constant and we also used the property $X_{u+h} - X_{u} \overset{d}{=} X_{h}$ for Brownian motion. As $k$ ranges within the inner sum in the second term in Eq. (\ref{lmpq0pf1}), there are $L$ values for which $(i-j)^{2} < (i-k)^{2}$, for some fixed $L$ in $\{1, \ldots, T \}$, and for the remaining $T-L-2$ values $(i-j)^{2} > (i-k)^{2}$. So, applying another fact mentioned in Remark 3, we get, \[ 2^{nd}\; \mathrm{term}\; \leq \frac{1}{(T-1)^{2}}\sum\limits_{\substack{j=1 \\ j \neq i}}^{T} \left\{\frac{L}{(i-j)^{2}} + C_{1}S_{2} \right\} \leq \frac{(C_{2} + C_{1}T)S_{2}}{(T-1)^{2}} = O(T^{-1}),  \]  where $C_{1}$ and $C_{2}$ are finite constants. Putting this together, we get Var($\hat{\gamma}_{i}^{0,0}$) $\rightarrow 0$ as $T \rightarrow \infty$. Using similar arguments it is easy to prove the remaining statements in the lemma. $\blacksquare$\\

\begin{lemma} $\mathrm{Var}(\hat{\gamma}_{i}^{q,p}$) $\rightarrow 0$, $\mathrm{Var}([\hat{\sigma}_{i}^{x,q,p}]^{2})$ $\rightarrow 0$ and $\mathrm{Var}([\hat{\sigma}_{i}^{y,q,p}]^{2})$ $\rightarrow 0$, for each $i$ and for each $q$ and $p$ with $0 < q \leq 1/2$, $p > 1/2$, as $T \rightarrow \infty$. \end{lemma}

\textbf{Proof:} \begin{align} \mathrm{Var}(\hat{\gamma}_{i}^{q,p}) &\leq \frac{1}{(T-1)^{2}}\sum\limits_{\substack{j=1 \\ j \neq i}}^{T}\frac{\mathbb{E}[(j^{q}X_{i}-X_{j}j^{-p})(j^{q}Y_{i}-Y_{j}j^{-p})]^{2}}{(i-j)^{4}}  \label{lmpqgnpf}  \\ &+ \frac{1}{(T-1)^{2}}\sum\limits_{\substack{j=1 \\ j \neq i}}^{T}\sum\limits_{\substack{k=1 \\ k \neq i,j}}^{T}\frac{\mathbb{E}[(j^{q}X_{i}-X_{j}j^{-p})(j^{q}Y_{i}-Y_{j}j^{-p})(k^{q}X_{i}-X_{k}k^{-p})(k^{q}Y_{i}-Y_{k}k^{-p})]}{(i-j)^{2}(i-k)^{2}}  \nonumber \end{align}
\begin{align*} [(j^{q}X_{i}-X_{j}j^{-p})(j^{q}Y_{i}-Y_{j}j^{-p})]^{2} &= [j^{2q}(X_{i}-X_{j})^{2} + 2X_{j}(X_{i}-X_{j})j^{q}(j^{q}-j^{-p}) + X_{j}^{2}(j^{q}-j^{-p})^{2}] \\\;\;\;\; &\times [j^{2q}(Y_{i}-Y_{j})^{2} + 2Y_{j}(Y_{i}-Y_{j})j^{q}(j^{q}-j^{-p}) + Y_{j}^{2}(j^{q}-j^{-p})^{2}] \end{align*} Expanding and taking expectations, many terms become zero by applying properties of Brownian motion. With the remaining non-zero terms, the first term in Eq. (\ref{lmpqgnpf}) becomes \begin{align} \frac{1}{(T-1)^{2}}\sum\limits_{\substack{j=1 \\ j \neq i}}^{T}(j-i)^{-4}&\left\{\mathbb{E}[j^{4q}(X_{j}-X_{i})^{2}(Y_{j}-Y_{i})^{2}] \right. \label{lmpqgnpf1term1} \\  \;\;\;&+\; \mathbb{E}[j^{2q}(X_{j}-X_{i})^{2}]\mathbb{E}[Y_{j}^{2}(j^{q}-j^{-p})^{2}] \label{lmpqgnpf1term2}\\  \;\;\;&+\; 4\mathbb{E}[X_{j}Y_{j}]\mathbb{E}[(X_{j}-X_{i})(Y_{j}-Y_{i})j^{2q}(j^{q}-j^{-p})^{2}] \label{lmpqgnpf1term3}\\  \;\;\;&+\;  \mathbb{E}[X_{j}^{2}]\mathbb{E}[(Y_{j}-Y_{i})^{2}j^{2q}(j^{q}-j^{-p})^{2}] \label{lmpqgnpf1term4}\\  \;\;\;&+\; \left. \mathbb{E}[X_{j}^{2}Y_{j}^{2}(j^{q}-j^{-p})^{4}]\right\}. \label{lmpqgnpf1term5} \end{align} We will now show that the product of $(j-i)^{-4}$ and each expectation term in Eq. (\ref{lmpqgnpf1term1}) to Eq. (\ref{lmpqgnpf1term5}) can be written as $C_{1} + C_{2}T$, for some finite constants $C_{1}$ and $C_{2}$ so that this product when multiplied by $(T-1)^{-2}$ is $O(T^{-1})$. Using Remark 3, for $j > i$, \begin{align} \frac{j^{4q}\mathbb{E}[(X_{j}-X_{i})^{2}(Y_{j}-Y_{i})^{2}]}{(j-i)^{4}} & \leq \frac{3j^{4q}}{(j-i)^{2}} < \frac{3j^{2}}{(j-i)^{2}}\; \mathrm{if}\; 4q \leq 2; \; \mathrm{i.e.\; if}\; q \leq 1/2 \nonumber \\ & \leq 3\left(1 + \frac{2i}{(j-i)} + \frac{i^{2}}{(j-i)^{2}} \right) < 12,\;\mathrm{if}\; j > 2i. \label{inb1} \end{align} Since the first $2i$ terms in the sum in (\ref{lmpqgnpf1term1}) is a constant, it is easy to see based on (\ref{inb1}) that the whole sum can be bounded by $C_{1} + C_{2}T$, for some finite constants $C_{1}$ and $C_{2}$. \[ \frac{j^{2q}(j^{q} - j^{-p})^{2}\mathbb{E}(X_{j} - X_{i})^{2}\mathbb{E}(Y_{j}^{2})}{(j-i)^{4}} < \frac{j^{2q+1}(j^{2q}+1)}{(j-i)^{3}} \] can be bounded by a constant if $4q + 1 \leq 3$; i.e. if $q \leq 1/2$. Hence the sum corresponding to (\ref{lmpqgnpf1term2}) can be expressed as $C_{1} + C_{2}T$. Similar reasoning can be applied to sums in (\ref{lmpqgnpf1term3}), (\ref{lmpqgnpf1term4}) and (\ref{lmpqgnpf1term5}) so that the first term in Eq. (\ref{lmpqgnpf}) is $O(T^{-1})$.

Expanding the second term in Eq. (\ref{lmpqgnpf}) gives 16 terms. The double-sum ranging over $j$ and $k$ for each of these 16 terms can be expressed as $[C_{1}S_{a} + C_{2}]S_{b}T$ for some $a > 1$ and $b > 1$, so that after multiplication with $(T-1)^{-2}$ each term becomes $O(T^{-1})$. For example, the last term will be \begin{align} \frac{\mathbb{E}(X_{j}Y_{j}j^{-2p}X_{k}Y_{k}k^{-2p})}{(j-i)^{2}(k-i)^{2}} & \leq \frac{j^{-2p+2}k^{-2p+2}}{(j-i)^{2}(k-i)^{2}}\times \frac{3}{j^{2}},\; \mathrm{for}\; k > j > i,  \nonumber \\ \;\mathrm{After\;summing\;over\;} k & < \frac{j^{-2p}}{(j-i)^{2}}(C + S_{2p}), \nonumber \\ \;\mathrm{After\;further\;summing\;over\;} j & < [C_{2} + C_{1}S_{2p}]S_{2p+2}T.  \nonumber \end{align}
Proceeding similarly with other terms, the 2$^{nd}$ term in Eq. (\ref{lmpqgnpf}) is also seen to be $O(T^{-1})$ and hence Var($\hat{\gamma}_{i}^{q,p}$) $\rightarrow 0$ as $T \rightarrow \infty$. Using similar arguments the remaining statements in the lemma can also be proved. $\blacksquare$\\

\begin{lemma} For each $t \in \{1, \ldots, n (= T)\}$, $ Q_{t} = \mathbb{E}(\hat{\gamma}_{t}^{q,p})/\sqrt{\mathbb{E}([\hat{\sigma}_{t}^{x,q,p}]^{2})
\mathbb{E}([\hat{\sigma}_{t}^{y,q,p}]^{2}) } \rightarrow \rho_{t}$ \textit{as} $T \rightarrow \infty, $ \textit{if} $p > q \geq 1/2$. \end{lemma}

\textbf{Proof:} \begin{align} (T-1)\mathbb{E}(\hat{\gamma}_{t}^{q,p}) & =  t\rho_{t}\sum\limits_{\substack{s=1 \\ s \neq t}}^{T}\frac{s^{2q}}{(s-t)^{2}}  - 2\left[\sum_{s=1}^{t-1} \frac{\rho_{s}s^{-p + q + 1}}{(t-s)^{2}} + \sum_{s=t+1}^{T}\frac{t\rho_{t}s^{-p + q}}{(s-t)^{2}} \right] + \sum\limits_{\substack{s=1 \\ s \neq t}}^{T}\frac{\rho_{s}s^{-2p+1}}{(s-t)^{2}} \nonumber \\  &= t\rho_{t}\sum\limits_{\substack{s=1 \\ s \neq t}}^{T}\frac{s^{2q}}{(s-t)^{2}} + \sum\limits_{\substack{s=1 \\ s \neq t}}^{T}\frac{\rho_{s}s^{-2p+1}}{(s-t)^{2}} - 2 \sum\limits_{\substack{s=1 \\ s \neq t}}^{T}\frac{\rho_{s}s^{-p + q + 1}}{(s-t)^{2}} + 2 \sum_{s=t+1}^{T} \frac{s^{-p + q}(s\rho_{s} - t\rho_{t})}{(s-t)^{2}}. \nonumber \end{align} Similarly, \begin{align} (T-1)\mathbb{E}([\hat{\sigma}_{t}^{x,p}]^{2}) = (T-1)\mathbb{E}([\hat{\sigma}_{t}^{y,p}]^{2}) &= t\sum\limits_{\substack{s=1 \\ s \neq t}}^{T}\frac{s^{2q}}{(s-t)^{2}} + \sum\limits_{\substack{s=1 \\ s \neq t}}^{T}\frac{s^{-2p+1}}{(s-t)^{2}} \nonumber \\ & - 2 \sum\limits_{\substack{s=1 \\ s \neq t}}^{T}\frac{s^{-p + q + 1}}{(s-t)^{2}} + 2 \sum_{s=t+1}^{T} \frac{s^{-p + q}(s - t)}{(s-t)^{2}}. \nonumber \end{align} Hence $Q_{t}$ can be written as a ratio \begin{equation} \label{Qt} Q_{t} = \frac{t\rho_{t} + K_{T}^{-1}A_{1,T} - 2K_{T}^{-1}A_{2,T} + 2K_{T}^{-1}A_{3,T} }{t + K_{T}^{-1}B_{1,T} - 2K_{T}^{-1}B_{2,T}  + 2K_{T}^{-1}B_{3,T}} ,  \end{equation} where
\[ K_{T} = \sum\limits_{\substack{s=1 \\ s \neq t}}^{T}\frac{s^{2q}}{(s-t)^{2}},\; A_{1,T} = \sum\limits_{\substack{s=1 \\ s \neq t}}^{T}\frac{\rho_{s}s^{-2p+1}}{(s-t)^{2}}, \; A_{2,T} = \sum\limits_{\substack{s=1 \\ s \neq t}}^{T}\frac{\rho_{s}s^{-p + q + 1}}{(s-t)^{2}},\; A_{3,T} = \sum_{s=t+1}^{T} \frac{s^{-p + q}(s\rho_{s} - t\rho_{t})}{(s-t)^{2}}, \] \[B_{1,T} = \sum\limits_{\substack{s=1 \\ s \neq t}}^{T}\frac{s^{-2p+1}}{(s-t)^{2}},\; B_{2,T} = \sum\limits_{\substack{s=1 \\ s \neq t}}^{T}\frac{s^{-p + q + 1}}{(s-t)^{2}},\; B_{3,T} = \sum_{s=t+1}^{T} \frac{s^{-p + q}(s - t)}{(s-t)^{2}}. \] $K_{T}$ diverges as $T \rightarrow \infty$ for $2q \geq 1$, and all other sums ($A_{1,T}$, $A_{2,T}$, $A_{3,T}$, $B_{1,T}$, $B_{2,T}$ and $B_{3,T}$) converges to finite constants in passage to the limit when $p > q \geq 1/2$. Thus by Eq. (\ref{Qt}), $Q_{t} \rightarrow \rho_{t}$ as $T \rightarrow \infty$ as desired. $\blacksquare$\\

\textbf{Proof of Theorem 3.1:} Note that the only value of $q$ for which \textit{both} Lemma 3.3 and Lemma 3.4 hold, is $q = 1/2$. Using Lemma 2.3, Chebychev's inequality, Slutsky's lemma, Continuous mapping theorem [8] and part (iii) of Theorem 2.7 in [8, p.10] we get $\hat{\rho}_{u}^{\,q,p} \rightarrow \mathbb{E}(\hat{\gamma}_{u}^{q,p})/\sqrt{\mathbb{E}([\hat{\sigma}_{u}^{x,q,p}]^{2})
\mathbb{E}([\hat{\sigma}_{u}^{y,q,p}]^{2}) }$ in probability for each $u$ and $p > q = 1/2$; this fact combined with Lemma 3.4 proves the theorem.   $\blacksquare$

\subsection{Weak consistency of dynamic correlation estimators in the Geometric Brownian Motion case}

The main results in this subsection are the following two theorems. The strategy for proofs in this subsection is the same as that used in the Brownian motion case; only the details in the calculations differ.\\

\begin{theorem} $\hat{\rho}_{t}^{a,b,c} \rightarrow \rho_{t}$ in probability for each $t$, $c > a > 0$ and $b > a + 10$, as $T \rightarrow \infty$, where $\hat{\rho}_{t}^{a,b,c} = \hat{\gamma}_{t}/(\hat{\sigma}_{t,W}\hat{\sigma}_{t,U})$ with $\hat{\gamma}_{t}$ given in Eq. (2) and $\hat{\sigma}_{t,W}$, $\hat{\sigma}_{t,U}$ the corresponding estimates of the variances. \end{theorem}

\begin{theorem} $\hat{\rho}_{t}^{a,b,c} \rightarrow \rho_{t}$ in probability for each $t$, $b > 15$ and $c > a > 0$, as $T \rightarrow \infty$, where $\hat{\rho}_{t}^{a,b,c} = \hat{\gamma}_{t}/(\hat{\sigma}_{t,W}\hat{\sigma}_{t,U})$ with $\hat{\gamma}_{t}$ given in Eq. (3) and $\hat{\sigma}_{t,W}$, $\hat{\sigma}_{t,U}$ the corresponding estimates of the variances. \end{theorem}

We first consider Theorem 3.5 and begin with a lemma similar to Lemma 3.3.\\

\begin{lemma} For $\hat{\gamma}_{t}$ given in Eq. (2), and the corresponding estimates of the variances $\hat{\sigma}_{t,W}^{2}$ and $\hat{\sigma}_{t,U}^{2}$, $\mathrm{Var}(\hat{\gamma}_{t}$) $\rightarrow 0$, $\mathrm{Var}([\hat{\sigma}_{t,W}]^{2})$ $\rightarrow 0$ and $\mathrm{Var}([\hat{\sigma}_{t,U}]^{2})$ $\rightarrow 0$, for each $t$ and for $c > a > 0$ and $b > a + 10$, as $T \rightarrow \infty$. \end{lemma}

Based on Eq. (2)
\begin{align} \displaystyle \hat{\gamma}_{t} = \frac{1}{e^{c\sigma^{2}T}}\sum_{k = 1}^{T}& \left\{ \left[e^{\frac{-b\sigma^{2}k}{2}}\left(e^{\sigma W_{k}}-e^{\frac{\sigma^{2}k}{2}}\right) - e^{\frac{a\sigma^{2}k}{2}}\left(e^{\sigma W_{t}}-e^{\frac{\sigma^{2}t}{2}}\right) \right] \right. \nonumber \\
& \left. \times \left[e^{\frac{-b\sigma^{2}k}{2}}\left(e^{\sigma U_{k}}-e^{\frac{\sigma^{2}k}{2}}\right) - e^{\frac{a\sigma^{2}k}{2}}\left(e^{\sigma U_{t}}-e^{\frac{\sigma^{2}t}{2}}\right) \right] \right\} \nonumber \\
& = G_{1} - G_{2} - G_{3} + G_{4}, \nonumber \end{align}
where
\begin{align} G_{1} & = e^{-c\sigma^{2}T}\sum_{k = 1}^{T}\left[e^{-b\sigma^{2}k}\left(e^{\sigma W_{k}}-e^{\frac{\sigma^{2}k}{2}}\right)\left(e^{\sigma U_{k}}-e^{\frac{\sigma^{2}k}{2}}\right) \right], \nonumber \\
G_{2} & = e^{-c\sigma^{2}T}\sum_{k = 1}^{T}\left[e^{(a-b)\sigma^{2}k/2}\left(e^{\sigma W_{k}}-e^{\frac{\sigma^{2}k}{2}}\right)\left(e^{\sigma U_{t}}-e^{\frac{\sigma^{2}t}{2}}\right) \right], \nonumber \\
G_{3} & = e^{-c\sigma^{2}T}\sum_{k = 1}^{T}\left[e^{(a-b)\sigma^{2}k/2}\left(e^{\sigma W_{t}}-e^{\frac{\sigma^{2}t}{2}}\right)\left(e^{\sigma U_{k}}-e^{\frac{\sigma^{2}k}{2}}\right) \right], \nonumber \\
G_{4} & = e^{-c\sigma^{2}T}\sum_{k = 1}^{T}\left[e^{a\sigma^{2}k}\left(e^{\sigma W_{t}}-e^{\frac{\sigma^{2}t}{2}}\right)\left(e^{\sigma U_{t}}-e^{\frac{\sigma^{2}t}{2}}\right) \right]. \nonumber
\end{align}

In order to show that $\mathrm{Var}(\hat{\gamma}_{t}) \rightarrow 0$, it suffices to show that $\mathbb{E} (\hat{\gamma}_{t}^{2}) \rightarrow 0$.

\begin{equation} \label{Gterms} \displaystyle \hat{\gamma}_{t}^{2} = [G_{1}^{2} + G_{2}^{2} + G_{3}^{2} + G_{4}^{2} -2G_{1}G_{2} - 2G_{1}G_{3} + 2G_{1}G_{4} + 2G_{2}G_{3} - 2G_{2}G_{4} - 2G_{3}G_{4}].  \end{equation} We will show that the expectation of each of the 10 terms in the right-hand side of Eq. (\ref{Gterms}) converges to zero as $T \rightarrow \infty$ for the range of $a, b$ and $c$ assumed in the statement of Lemma 3.7, thereby showing that $\mathbb{E}(\hat{\gamma}_{t}^{2}) \rightarrow 0, \;\mathrm{as}\; T \rightarrow \infty$. The proof is lengthy but the calculations involved are very routine and repetitious. Hence we relegate the proof of Lemma 3.7 to the appendix.

Next we state and prove a lemma similar to Lemma 3.4.
\begin{lemma} For $\hat{\gamma}_{t}$ given in Eq. (2), and the corresponding estimates of the variances $\hat{\sigma}_{t,W}^{2}$ and $\hat{\sigma}_{t,U}^{2}$,
 $\mathbb{E}(\hat{\gamma}_{t})/\sqrt{\mathbb{E}(\hat{\sigma}_{t,W}^{2})\mathbb{E}(\hat{\sigma}_{t,U}^{2})} \rightarrow \rho_{t}$, as $T \rightarrow \infty$, if $b > a > 0$. \end{lemma}

\textbf{Proof of Lemma 3.8}: $e^{c\sigma^{2}T}\hat{\gamma}_{t} = \sum_{k=1}^{t-1}\{\cdot\} + \sum_{k=t}^{T}\{\cdot\}$ where \begin{align} \{\cdot\} &= \left[e^{-b\sigma^{2}k/2}(e^{\sigma W_{k}} - e^{\sigma^{2}k/2}) - e^{a\sigma^{2}k/2}(e^{\sigma W_{t}} - e^{\sigma^{2}t/2}) \right]  \nonumber \\ & \;\;\;\;\;\;\;\;\;\; \times \left[e^{-b\sigma^{2}k/2}(e^{\sigma U_{k}} - e^{\sigma^{2}k/2}) - e^{a\sigma^{2}k/2}(e^{\sigma U_{t}} - e^{\sigma^{2}t/2}) \right]\nonumber \\
& = e^{-b\sigma^{2}k}\left[e^{\sigma(W_{k} + U_{k})} - e^{\sigma^{2}k/2}e^{\sigma W_{k}} - e^{\sigma^{2}k/2}e^{\sigma U_{k}} + e^{\sigma^{2}k} \right] \nonumber \\ & \;\;\;\;\;\;\;\;\;\;- e^{(a-b)\sigma^{2}k/2}\left[e^{\sigma(W_{k} + U_{t})} - e^{\sigma^{2}k/2}e^{\sigma U_{t}} - e^{\sigma^{2}t/2}e^{\sigma W_{k}} + e^{\sigma^{2}(k+t)/2} \right] \nonumber \\ & \;\;\;\;\;\;\;\;\;\;- e^{(a-b)\sigma^{2}k/2}\left[e^{\sigma(U_{k} + W_{t})} - e^{\sigma^{2}t/2}e^{\sigma U_{k}} - e^{\sigma^{2}k/2}e^{\sigma W_{t}} + e^{\sigma^{2}(k+t)/2} \right] \nonumber \\ & \;\;\;\;\;\;\;\;\;\;+  e^{a\sigma^{2}k}\left[e^{\sigma(W_{t} + U_{t})} - e^{\sigma^{2}t/2}e^{\sigma W_{t}} - e^{\sigma^{2}t/2}e^{\sigma U_{t}} + e^{\sigma^{2}t} \right] \nonumber \end{align} We recall $\mathbb{E}(e^{\sigma(W_{k} + U_{k})}) = e^{\sigma^{2}k}[1 + \rho_{k}(e^{\sigma^{2}k}-1)]$, $\mathbb{E}(e^{\sigma(W_{t} + U_{t})}) = e^{\sigma^{2}t}[1 + \rho_{t}(e^{\sigma^{2}t}-1)]$, $W_{k} = \rho_{k}U_{k} + \sqrt{1-\rho_{k}^{2}}M_{k}^{(1)}$ and $\mathbb{E}(e^{\sigma(W_{k} + U_{t})}) = \mathbb{E}\left(\exp\{\sigma[\rho_{k}U_{k} + U_{t} + \sqrt{1-\rho_{k}^{2}}M_{k}^{(1)}] \} \right)$ so that when $k \geq t$, \begin{align} \mathbb{E}(e^{\sigma(W_{k} + U_{t})}) & = \mathbb{E}(e^{\sigma\rho_{k}(U_{k}-U_{t})})\mathbb{E}(e^{\sigma(1+\rho_{k})U_{t}})\mathbb{E}(e^{\sigma\sqrt{1-\rho_{k}^{2}}M_{k}^{(1)}}) \nonumber \\ & = e^{ \frac{\sigma^{2}}{2}[\rho_{k}^{2}(k-t) + (1 + \rho_{k})^{2}t + (1-\rho_{k}^{2})k ] }  = e^{\sigma^{2}(k+t)/2}e^{\sigma^{2}\rho_{k}t} \leq e^{\sigma^{2}k/2}e^{3\sigma^{2}t/2}. \nonumber  \end{align} The same relationship as above holds for $\mathbb{E}(e^{\sigma(U_{k} + W_{t})})$ also. When $k \leq t$, \begin{align} \mathbb{E}(e^{\sigma(W_{k} + U_{t})}) & = \mathbb{E}(e^{\sigma\rho_{k}(U_{t}-U_{k})})\mathbb{E}(e^{\sigma(1+\rho_{k})U_{k}})\mathbb{E}(e^{\sigma\sqrt{1-\rho_{k}^{2}}M_{k}^{(1)}}) \nonumber \\ & = e^{ \frac{\sigma^{2}}{2}[\rho_{k}^{2}(t-k) + (1 + \rho_{k})^{2}t + (1-\rho_{k}^{2})k ] }  = e^{\sigma^{2}(k+t)/2}e^{\sigma^{2}\rho_{k}k} \leq e^{\sigma^{2}t/2}e^{3\sigma^{2}k/2}; \nonumber  \end{align} again the same relationship holding for $\mathbb{E}(e^{\sigma(U_{k} + W_{t})})$ also. Hence \begin{align} e^{c\sigma^{2}T}\mathbb{E}(\hat{\gamma}_{t}) &= \sum_{k=1}^{T}e^{-b\sigma^{2}k}\{e^{\sigma^{2}k}[1 + \rho_{k}(e^{\sigma^{2}k}-1)] - e^{\sigma^{2}k} \} + \sum_{k=1}^{T}e^{a\sigma^{2}k}\{e^{\sigma^{2}t}[1 + \rho_{t}(e^{\sigma^{2}t}-1)] - e^{\sigma^{2}t} \} \nonumber \\ &- 2 \sum_{k=1}^{t-1}e^{(a-b)\sigma^{2}k/2}\{e^{\sigma^{2}(k+t)/2}[e^{\sigma^{2}\rho_{k}k}-1]\} - 2 \sum_{k=t}^{T}e^{(a-b)\sigma^{2}k/2}\{e^{\sigma^{2}(k+t)/2}[e^{\sigma^{2}\rho_{k}t}-1]\} \nonumber \\ & = \left\{ \sum_{k=1}^{T}\rho_{k}[e^{(2-b)\sigma^{2}k} - e^{(1-b)\sigma^{2}k}] \right\} + \rho_{t}(e^{\sigma^{2}t}-1)e^{(a+t)\sigma^{2}}\left[\frac{e^{a\sigma^{2}T}-1}{e^{a\sigma^{2}}-1} \right] \nonumber \\ & -2\sum_{k=1}^{t-1}e^{(a-b)\sigma^{2}k/2}\{e^{\sigma^{2}(k+t)/2}[e^{\sigma^{2}\rho_{k}k}-e^{\sigma^{2}\rho_{k}t}]\} - 2 \sum_{k=1}^{T}e^{(a-b)\sigma^{2}k/2}\{e^{\sigma^{2}(k+t)/2}[e^{\sigma^{2}\rho_{k}t}-1]\} \nonumber \\ & = A_{T} + D_{T}\rho_{t} - 2B_{t} - 2C_{T}, \nonumber  \end{align} where $B_{t}$ is a finite sum of ($t-1$) terms that does not depend on $T$, $A_{T}$ and $C_{T}$ converge to finite constants as $T \rightarrow \infty$ since we assume $b > a$, and \[ D_{T} = e^{(a+t)\sigma^{2}}(e^{\sigma^{2}t}-1)\left[\frac{e^{a\sigma^{2}T}-1}{e^{a\sigma^{2}}-1} \right]\; \mathrm{diverges\; as}\; T \rightarrow \infty. \]

Similarly, $e^{c\sigma^{2}T}\hat{\sigma}_{t,W}^{2} = \sum_{k=1}^{t-1}\{\cdot\} + \sum_{k=t}^{T}\{\cdot\}$ where \begin{align} \{\cdot\} &= \left[e^{-b\sigma^{2}k/2}(e^{\sigma W_{k}} - e^{\sigma^{2}k/2}) - e^{a\sigma^{2}k/2}(e^{\sigma W_{t}} - e^{\sigma^{2}t/2}) \right]^{2} \nonumber \\ & = e^{-b\sigma^{2}k}(e^{\sigma W_{k}} - e^{\sigma^{2}k/2})^{2} + e^{a\sigma^{2}k}(e^{\sigma W_{t}} - e^{\sigma^{2}t/2})^{2} - 2e^{(a-b)\sigma^{2}k/2}(e^{\sigma W_{k}} - e^{\sigma^{2}k/2})(e^{\sigma W_{t}} - e^{\sigma^{2}t/2}). \nonumber \end{align} We have $\mathbb{E}(e^{\sigma W_{k}} - e^{\sigma^{2}k/2})^{2} = \mathbb{E}[e^{2\sigma W_{k}} - 2e^{\sigma W_{k}}e^{\sigma^{2}k/2} + e^{\sigma^{2}k}] = e^{2\sigma^{2}k} - e^{\sigma^{2}k} = e^{\sigma^{2}k}(e^{\sigma^{2}k}-1)$ and $\mathbb{E}(e^{\sigma W_{t}} - e^{\sigma^{2}t/2})^{2} = e^{\sigma^{2}t}(e^{\sigma^{2}t}-1)$.

When $k \geq t, \mathbb{E}(e^{\sigma (W_{k} + W_{t})})  = \mathbb{E}(e^{\sigma(W_{k}-W_{t})})\mathbb{E}(e^{2\sigma W_{t}}) = \exp \{\frac{\sigma^{2}}{2}(k-t) + 2\sigma^{2}t \} = e^{\sigma^{2}(k+t)/2}e^{\sigma^{2}t}$, and when $k \leq t$, $\mathbb{E}(e^{\sigma (W_{k} + W_{t})}) = e^{\sigma^{2}(k+t)/2}e^{\sigma^{2}k}.$ Putting it all together, \begin{align} e^{c\sigma^{2}T}\mathbb{E}(\hat{\sigma}_{t,W}^{2}) & = \sum_{k=1}^{T}[e^{(2-b)\sigma^{2}k} - e^{(1-b)\sigma^{2}k}] + e^{\sigma^{2}t}(e^{\sigma^{2}t}-1)e^{a\sigma^{2}}\left[\frac{e^{a\sigma^{2}T}-1}{e^{a\sigma^{2}}-1} \right] \nonumber \\ &\;\;\;\; -2\sum_{k=1}^{t-1}e^{(a-b)\sigma^{2}k/2}e^{\sigma^{2}(k+t)/2}[e^{\sigma^{2}k}-e^{\sigma^{2}t}] - 2e^{\sigma^{2}t/2}(e^{\sigma^{2}t}-1)\sum_{k=1}^{T}e^{\sigma^{2}k(a-b+1)/2}. \nonumber \end{align} Hence $e^{c\sigma^{2}T}\mathbb{E}(\hat{\sigma}_{t,W}^{2})$ can be written as $E_{T} + D_{T} - 2F_{t} - 2G_{T}$, where $D_{T}$ is same as the one given further above, $F_{t}$ is a finite sum that does not depend on $T$ and $E_{T}$ and $G_{T}$ converges to finite constants as $T \rightarrow \infty$. The same result holds for $\mathbb{E}(\hat{\sigma}_{t,U}^{2})$ also. Hence \begin{align} \mathbb{E}(\hat{\gamma}_{t})/\sqrt{\mathbb{E}(\hat{\sigma}_{t,W}^{2})\mathbb{E}(\hat{\sigma}_{t,U}^{2})} &= \frac{\rho_{t} + A_{T}/D_{T} - 2B_{t}/D_{T} - 2C_{T}/D_{T}}{1 + E_{T}/D_{T} - 2F_{t}/D_{T} - 2G_{T}/D_{T}} \nonumber \\ &\rightarrow \rho_{t},\; \mathrm{as}\; T \rightarrow \infty, \nonumber\end{align} since $D_{T}$ diverges as $T \rightarrow \infty$, $B_{t}$ and $F_{t}$ are finite terms that do not depend on $T$ and $A_{T}, C_{T}, E_{T}$ and $G_{T}$ are all finite constants in the passage to the limit. This proves the lemma. $\blacksquare$\\

Lemmas 3.7 and 3.8 suffices to prove Theorem 3.5. Similarly lemmas 3.9 and 3.10 given below will suffice to prove Theorem 3.6.\\

\begin{lemma} For $\hat{\gamma}_{t}$ given in Eq. (3), and the corresponding estimates of the variances $\hat{\sigma}_{t,W}^{2}$ and $\hat{\sigma}_{t,U}^{2}$, $\mathrm{Var}(\hat{\gamma}_{t}$) $\rightarrow 0$, $\mathrm{Var}([\hat{\sigma}_{t,W}]^{2})$ $\rightarrow 0$ and $\mathrm{Var}([\hat{\sigma}_{t,U}]^{2})$ $\rightarrow 0$, for each $t$ and for $c > a > 0$ and $b > 15$, as $T \rightarrow \infty$. \end{lemma}

\textbf{Proof of Lemma 3.9:}

\begin{align}
\hat{\gamma}_{t}^{2} = e^{-2c\sigma^{2}T}&\left\{
\sum_{k=1}^{T}\sum_{j=1}^{T} e^{a\sigma^{2}(k+j)}(e^{\sigma W_{t}} - e^{\sigma^{2}t/2})^{2}(e^{\sigma U_{t}} - e^{\sigma^{2}t/2})^{2} \right. \nonumber \\
&-\sum_{k=1}^{T}\sum_{j=1}^{T} e^{a\sigma^{2}k}e^{-b\sigma^{2}j}(e^{\sigma W_{t}} - e^{\sigma^{2}t/2})(e^{\sigma U_{t}} - e^{\sigma^{2}t/2})(e^{\sigma W_{j}} - e^{\sigma^{2}j/2})(e^{\sigma U_{j}} - e^{\sigma^{2}j/2})  \nonumber \\ &-\sum_{k=1}^{T}\sum_{j=1}^{T} e^{a\sigma^{2}k}e^{-b\sigma^{2}j}(e^{\sigma W_{k}} - e^{\sigma^{2}k/2})(e^{\sigma U_{k}} - e^{\sigma^{2}k/2})(e^{\sigma W_{t}} - e^{\sigma^{2}t/2})(e^{\sigma U_{t}} - e^{\sigma^{2}t/2})  \nonumber \\ &+ \left. \sum_{k=1}^{T}\sum_{j=1}^{T} e^{-b\sigma^{2}(k+j)}(e^{\sigma W_{k}} - e^{\sigma^{2}k/2})(e^{\sigma U_{k}} - e^{\sigma^{2}k/2})(e^{\sigma W_{j}} - e^{\sigma^{2}j/2})(e^{\sigma U_{j}} - e^{\sigma^{2}j/2}) \right\} \nonumber
\end{align}

It is easy to see that when taking expectations, the first term converges to zero as $T \rightarrow \infty$ if $c > a$. The expectations of the second and third terms are the same. We exhibit only the calculation for the third term. Within the third term if we expand \[ (e^{\sigma W_{k}} - e^{\sigma^{2}k/2})(e^{\sigma U_{k}} - e^{\sigma^{2}k/2})(e^{\sigma W_{t}} - e^{\sigma^{2}t/2})(e^{\sigma U_{t}} - e^{\sigma^{2}t/2}) \] there will be 16 terms. We exhibit only the calculation for the term containing $e^{\sigma(W_{k} + U_{k} + W_{t} + U_{t})}$ since similar calculations apply to the remaining 15 terms.

Recalling again $W_{k} = \rho_{k}U_{k}  + \sqrt{1-\rho_{k}^{2}}M_{k}^{(1)}$ we have \begin{align} W_{k} + U_{k} + W_{t} + U_{t} &\leq (1 + \rho_{k})U_{k} + \sqrt{1-\rho_{k}^{2}}M_{k}^{(1)} + (1 + \rho_{t})U_{t} + \sqrt{1-\rho_{t}^{2}}M_{t}^{(1)} \nonumber \\ &\leq 2|U_{k}| + |M_{k}^{(1)}| + 2|U_{t}| + |M_{t}^{(1)}| \nonumber \end{align} so that when $k \geq t$, \begin{align}W_{k} + U_{k} + W_{t} + U_{t} &\leq 2|U_{k} - U_{t}| + 4|U_{t}| + |M_{k}^{(1)} - M_{t}^{(1)}| + 2|M_{t}^{(1)}| \nonumber \\ \mathbb{E}(e^{\sigma(W_{k} + U_{k} + W_{t} + U_{t})}) &\leq C\exp\{\frac{\sigma^{2}}{2}[4(k-t) + 16t + (k-t) + 4t ] \} \leq Ce^{5\sigma^{2}k} \nonumber  \end{align} and when $k \leq t$, \begin{align}W_{k} + U_{k} + W_{t} + U_{t} &\leq 2|U_{t} - U_{k}| + 4|U_{k}| + |M_{t}^{(1)} - M_{k}^{(1)}| + 2|M_{k}^{(1)}| \nonumber \\ \mathbb{E}(e^{\sigma(W_{k} + U_{k} + W_{t} + U_{t})}) &\leq C\exp\{\frac{\sigma^{2}}{2}[4(t-k) + 16k + (t-k) + 4k ] \} \leq Ce^{15\sigma^{2}k}. \nonumber  \end{align} Hence the expectation of the term containing $e^{\sigma(W_{k} + U_{k} + W_{t} + U_{t})}$ in the third term is \begin{align}
 &\leq C e^{-2c\sigma^{2}T}  \sum_{k=1}^{T}\sum_{j=1}^{T} e^{a\sigma^{2}k}e^{-b\sigma^{2}j}e^{15\sigma^{2}k}  \nonumber \\ &= C e^{-2c\sigma^{2}T}  e^{a\sigma^{2}}\left[\frac{e^{a\sigma^{2}T}-1}{e^{a\sigma^{2}}-1} \right]e^{(15-b)\sigma^{2}}\left[\frac{1 - e^{(15-b)\sigma^{2}T}}{1 - e^{(15-b)\sigma^{2}}} \right]  \nonumber \\ &\leq C[e^{(a-2c)\sigma^{2}T} - e^{-2c\sigma^{2}T}][1 - e^{(15-b)\sigma^{2}T}] \nonumber \\ & \rightarrow 0, \;\mathrm{as}\; T \rightarrow \infty,\;\mathrm{if}\; c > a/2 > 0\; \mathrm{and}\; b > 15. \nonumber
\end{align}

Finally we exhibit the calculation for the first term in the expansion of the  fourth term of $\hat{\gamma}_{t}^{2}$; the corresponding expectation is \begin{align} &= e^{-2c\sigma^{2}T}\sum_{k=1}^{T}\sum_{j=1}^{T} e^{-b\sigma^{2}(k+j)}\mathbb{E}(e^{\sigma (W_{k} + U_{k} + W_{j} + U_{j})} ) \nonumber \\ & \leq C e^{-2c\sigma^{2}T}\sum_{k=1}^{T}\sum_{j=1}^{T}e^{(15-b)\sigma^{2}(k+j)} \nonumber \\ & \rightarrow 0, \;\mathrm{as}\; T \rightarrow \infty,\;\mathrm{if}\; c > 0\; \mathrm{and}\; b > 15. \nonumber
\end{align} The calculations for the remaining terms within the expansion of the fourth term are similar. $\blacksquare$\\

\begin{lemma} For $\hat{\gamma}_{t}$ given in Eq. (3), and the corresponding estimates of the variances $\hat{\sigma}_{t,W}^{2}$ and $\hat{\sigma}_{t,U}^{2}$,
 $\mathbb{E}(\hat{\gamma}_{t})/\sqrt{\mathbb{E}(\hat{\sigma}_{t,W}^{2})\mathbb{E}(\hat{\sigma}_{t,U}^{2})} \rightarrow \rho_{t}$, as $T \rightarrow \infty$, if $b > 2$ and $a > 0$. \end{lemma}

\textbf{Proof of Lemma 3.10}:
\[\displaystyle e^{c\sigma^{2}T}\hat{\gamma}_{t} = \sum_{k=1}^{T} \left\{e^{a\sigma^{2}k}[(e^{\sigma W_{t}} - e^{\sigma^{2}t/2})(e^{\sigma U_{t}} - e^{\sigma^{2}t/2}) ] - e^{-b\sigma^{2}k}[(e^{\sigma W_{k}} - e^{\sigma^{2}k/2})(e^{\sigma U_{k}} - e^{\sigma^{2}k/2})  ] \right\} \]
\begin{align} e^{c\sigma^{2}T}\mathbb{E}(\hat{\gamma}_{t}) &= \sum_{k=1}^{T} \left\{ e^{a\sigma^{2}k}\left[e^{\sigma^{2}t}[1 + \rho_{t}(e^{\sigma^{2}t}-1)] - e^{\sigma^{2}t}\right] \right. \nonumber \\ & \left. \;\;\;\;\;\;\;\;\;\;\;\;\;\;- e^{-b\sigma^{2}k}\left[e^{\sigma^{2}k}[1 + \rho_{k}(e^{\sigma^{2}k}-1)] - e^{\sigma^{2}k}\right] \right\} \nonumber \\ &= \rho_{t}e^{\sigma^{2}t}(e^{\sigma^{2}t} - 1)e^{a\sigma^{2}}\left[\frac{e^{a\sigma^{2}T}-1}{e^{a\sigma^{2}}-1} \right] - \sum_{k=1}^{T} \left\{ e^{(1-b)\sigma^{2}k}\rho_{k}(e^{\sigma^{2}k}-1) \right\} \nonumber \\ & = D_{T}\rho_{t} - A_{T}, \mathrm{where}\; D_{T} \; \mathrm{diverges}\; \mathrm{and}\; A_{T}\; \mathrm{converges}\; \mathrm{if}\; b > 2,\; \mathrm{as}\; T \rightarrow \infty. \nonumber \end{align}

\[\displaystyle e^{c\sigma^{2}T}\hat{\sigma}_{t,W}^{2} = \sum_{k=1}^{T} \left\{e^{a\sigma^{2}k}[(e^{\sigma W_{t}} - e^{\sigma^{2}t/2})]^{2} - e^{-b\sigma^{2}k}[(e^{\sigma W_{k}} - e^{\sigma^{2}k/2}) ]^{2} \right\} \]
\begin{align} e^{c\sigma^{2}T}\mathbb{E}(\hat{\sigma}_{t,W}^{2}) & = \sum_{k=1}^{T} \left\{e^{a\sigma^{2}k}e^{\sigma^{2}t}(e^{\sigma^{2}t}-1) - e^{(1-b)\sigma^{2}k}(e^{\sigma^{2}k}-1) \right\} \nonumber \\ & = D_{T}\rho_{t} - B_{T}, \; \mathrm{where}\; B_{T}\; \mathrm{converges}\; \mathrm{if}\; b > 2. \nonumber \end{align} The same result holds for $e^{c\sigma^{2}T}\mathbb{E}(\hat{\sigma}_{t,W}^{2})$ also. Hence \begin{align}\mathbb{E}(\hat{\gamma}_{t})/\sqrt{\mathbb{E}(\hat{\sigma}_{t,W}^{2})\mathbb{E}(\hat{\sigma}_{t,U}^{2})} & = \frac{D_{T}\rho_{t} - A_{T}}{D_{T} - B_{T}} = \frac{\rho_{t} - (A_{T}/D_{T})}{1 - (B_{T}/D_{T})} \nonumber \\ &\rightarrow \rho_{t},\; \mathrm{as}\; T \rightarrow \infty,\; \mathrm{if}\; a > 0\; \mathrm{and}\; b > 2. \nonumber\end{align} This proves the lemma. $\blacksquare$\\

\section{Conclusions}

In this paper we presented estimators for dynamic correlation between a pair of correlated Brownian motions and separately estimators for dynamic correlation between a pair of correlated Geometric Brownian motions. The main thrust of this paper was in showing the weak consistency of all the estimators presented.

The type of estimators that we presented here are generalizations of estimators that we introduced in an earlier work [4]. The special case that we considered in [4] had shown good empirical performance in a few simulation scenarios. Assessment of empirical performance of generalized versions introduced in this paper will be of interest for future work. The `hyperparameters' $p, q, a, b$ and $c$ occurring in the estimators presented in this paper may be considered as tuning parameters, similar to considering the `window-size' as a tuning parameter in the dynamic correlation estimation approach based on sliding windows. Computational approaches such as cross-validation may be utilized in obtaining the optimal hyperparameters. The empirical study that we will conduct in future will help to provide a better understanding on those topics as well.

\section{Appendix}

We provide the proof of Lemma 3.7 for the convenience of the reader.

\textbf{Proof of Lemma 3.7:}

Based on Eq. (2)
\begin{align} \displaystyle \hat{\gamma}_{t} = \frac{1}{e^{c\sigma^{2}T}}\sum_{k = 1}^{T}& \left\{ \left[e^{\frac{-b\sigma^{2}k}{2}}\left(e^{\sigma W_{k}}-e^{\frac{\sigma^{2}k}{2}}\right) - e^{\frac{a\sigma^{2}k}{2}}\left(e^{\sigma W_{t}}-e^{\frac{\sigma^{2}t}{2}}\right) \right] \right. \nonumber \\
& \left. \times \left[e^{\frac{-b\sigma^{2}k}{2}}\left(e^{\sigma U_{k}}-e^{\frac{\sigma^{2}k}{2}}\right) - e^{\frac{a\sigma^{2}k}{2}}\left(e^{\sigma U_{t}}-e^{\frac{\sigma^{2}t}{2}}\right) \right] \right\} \nonumber \\
& = G_{1} - G_{2} - G_{3} + G_{4}, \nonumber \end{align}
where
\begin{align} G_{1} & = e^{-c\sigma^{2}T}\sum_{k = 1}^{T}\left[e^{-b\sigma^{2}k}\left(e^{\sigma W_{k}}-e^{\frac{\sigma^{2}k}{2}}\right)\left(e^{\sigma U_{k}}-e^{\frac{\sigma^{2}k}{2}}\right) \right], \nonumber \\
G_{2} & = e^{-c\sigma^{2}T}\sum_{k = 1}^{T}\left[e^{(a-b)\sigma^{2}k/2}\left(e^{\sigma W_{k}}-e^{\frac{\sigma^{2}k}{2}}\right)\left(e^{\sigma U_{t}}-e^{\frac{\sigma^{2}t}{2}}\right) \right], \nonumber \\
G_{3} & = e^{-c\sigma^{2}T}\sum_{k = 1}^{T}\left[e^{(a-b)\sigma^{2}k/2}\left(e^{\sigma W_{t}}-e^{\frac{\sigma^{2}t}{2}}\right)\left(e^{\sigma U_{k}}-e^{\frac{\sigma^{2}k}{2}}\right) \right], \nonumber \\
G_{4} & = e^{-c\sigma^{2}T}\sum_{k = 1}^{T}\left[e^{a\sigma^{2}k}\left(e^{\sigma W_{t}}-e^{\frac{\sigma^{2}t}{2}}\right)\left(e^{\sigma U_{t}}-e^{\frac{\sigma^{2}t}{2}}\right) \right]. \nonumber
\end{align}

In order to show that $\mathrm{Var}(\hat{\gamma}_{t}) \rightarrow 0$, it suffices to show that $\mathbb{E} (\hat{\gamma}_{t}^{2}) \rightarrow 0$.

\begin{equation} \label{G.gamma} \displaystyle \hat{\gamma}_{t}^{2} = [G_{1}^{2} + G_{2}^{2} + G_{3}^{2} + G_{4}^{2} -2G_{1}G_{2} - 2G_{1}G_{3} + 2G_{1}G_{4} + 2G_{2}G_{3} - 2G_{2}G_{4} - 2G_{3}G_{4}]. \end{equation} We will show that the expectation of each of the 10 terms in the right-hand side of (\ref{G.gamma}) converges to zero as $T \rightarrow \infty$ for the range of $a, b$ and $c$ assumed in the statement of Lemma 3.7, thereby showing that $\mathbb{E}(\hat{\gamma}_{t}^{2}) \rightarrow 0, \;\mathrm{as}\; T \rightarrow \infty$.

\begin{align}\displaystyle  G_{1}^{2} = e^{-2c\sigma^{2}T}& \left\{\sum_{k=1}^{T}\sum_{j=1}^{T}e^{-b\sigma^{2}(k+j)}e^{\sigma(W_{k} + U_{k} + W_{j} + U_{j})} - \sum_{k=1}^{T}\sum_{j=1}^{T}e^{-b\sigma^{2}(k+j)}e^{\frac{\sigma^{2}j}{2}}e^{\sigma(W_{k} + U_{k} + W_{j})}  \right. \nonumber \\
& - \sum_{k=1}^{T}\sum_{j=1}^{T}e^{-b\sigma^{2}(k+j)}e^{\frac{\sigma^{2}j}{2}}e^{\sigma(W_{k} + U_{k} + U_{j})}  +  \sum_{k=1}^{T}\sum_{j=1}^{T}e^{-b\sigma^{2}(k+j)}e^{\frac{\sigma^{2}j}{2}}e^{\sigma(W_{k} + U_{k})} \nonumber  \\
& - \sum_{k=1}^{T}\sum_{j=1}^{T}e^{-b\sigma^{2}(k+j)}e^{\frac{\sigma^{2}k}{2}}e^{\sigma(W_{k} + W_{j} + U_{j})}  +  \sum_{k=1}^{T}\sum_{j=1}^{T}e^{-b\sigma^{2}(k+j)}e^{\frac{\sigma^{2}(k+j)}{2}}e^{\sigma(W_{k} + W_{j})} \nonumber \\
& + \sum_{k=1}^{T}\sum_{j=1}^{T}e^{-b\sigma^{2}(k+j)}e^{\frac{\sigma^{2}(k+j)}{2}}e^{\sigma(W_{k} + U_{j})}  - \sum_{k=1}^{T}\sum_{j=1}^{T}e^{-b\sigma^{2}(k+j)}e^{\frac{\sigma^{2}(k+2j)}{2}}e^{\sigma W_{k}} \nonumber \\
& - \sum_{k=1}^{T}\sum_{j=1}^{T}e^{-b\sigma^{2}(k+j)}e^{\frac{\sigma^{2}k}{2}}e^{\sigma(U_{k} + W_{j} + U_{j})} + \sum_{k=1}^{T}\sum_{j=1}^{T}e^{-b\sigma^{2}(k+j)}e^{\frac{\sigma^{2}(k+j)}{2}}e^{\sigma( U_{k} + W_{j})} \nonumber \\
& + \sum_{k=1}^{T}\sum_{j=1}^{T}e^{-b\sigma^{2}(k+j)}e^{\frac{\sigma^{2}(k+j)}{2}}e^{\sigma( U_{k} + U_{j})}  - \sum_{k=1}^{T}\sum_{j=1}^{T}e^{-b\sigma^{2}(k+j)}e^{\frac{\sigma^{2}(k+2j)}{2}}e^{\sigma U_{k}} \nonumber \\
& + \sum_{k=1}^{T}\sum_{j=1}^{T}e^{-b\sigma^{2}(k+j)}e^{\sigma^{2}k}e^{\sigma( W_{j} + U_{j})}  - \sum_{k=1}^{T}\sum_{j=1}^{T}e^{-b\sigma^{2}(k+j)}e^{\frac{\sigma^{2}(2k+j)}{2}}e^{\sigma W_{j}} \nonumber \\
& \left. - \sum_{k=1}^{T}\sum_{j=1}^{T}e^{-b\sigma^{2}(k+j)}e^{\frac{\sigma^{2}(2k+j)}{2}}e^{\sigma U_{j}} + \sum_{k=1}^{T}\sum_{j=1}^{T}e^{-b\sigma^{2}(k+j)}e^{\sigma^{2}(k+j)} \right\} \nonumber
\end{align}

We label the 16 terms in the right-hand side of the equation for $G_{1}^{2}$ as $G_{1}^{2}.\mathrm{term}1$, $G_{1}^{2}.\mathrm{term}2$, ..., $G_{1}^{2}.\mathrm{term}16$, in the same order as they appear above, and we will show the expectation of each of these terms to converge to zero as $T \rightarrow \infty$. We start with the last term since it is the easiest.

\begin{align}\displaystyle  \mathbb{E}(G_{1}^{2}.\mathrm{term16}) &= e^{-2c\sigma^{2}T}\sum_{k=1}^{T}\sum_{j=1}^{T}e^{-b\sigma^{2}(k+j)}e^{\sigma^{2}(k+j)} = e^{-2c\sigma^{2}T}e^{2(1-b)\sigma^{2}}\left[\frac{1 - e^{(1-b)\sigma^{2}T}}{1 - e^{(1-b)\sigma^{2}}} \right]^{2} \nonumber \\ & \rightarrow 0\;\mathrm{as}\; T \rightarrow \infty, \;\mathrm{if}\; c > 0\;\mathrm{and}\; b > 1. \nonumber \end{align}

Next we show that the expectation of $G_{1}^{2}.\mathrm{term}8$ converges to zero. The same proof holds for the terms $G_{1}^{2}.\mathrm{term}12, G_{1}^{2}.\mathrm{term}14$ and $G_{1}^{2}.\mathrm{term}15$ also.

\begin{align}\displaystyle  \mathbb{E}(G_{1}^{2}.\mathrm{term8}) &= e^{-2c\sigma^{2}T} \sum_{k=1}^{T}\sum_{j=1}^{T}e^{-b\sigma^{2}(k+j)}e^{\frac{\sigma^{2}(k+2j)}{2}}\mathbb{E}(e^{\sigma(W_{k})}) = e^{-2c\sigma^{2}T}\sum_{k=1}^{T}\sum_{j=1}^{T}e^{(1-b)\sigma^{2}(k+j)} \nonumber \\ &= e^{-2c\sigma^{2}T} e^{2(1-b)\sigma^{2}}\left[\frac{1 - e^{2(1-b)\sigma^{2}T}}{1 - e^{2(1-b)\sigma^{2}}} \right]^{2} \rightarrow 0\;\mathrm{as}\; T \rightarrow \infty, \;\mathrm{if}\; c > 0\;\mathrm{and}\; b > 1. \nonumber \end{align}

We will use the following inequality frequently in the remaining parts of the proof of Lemma 3.7.
\begin{equation} \label{absineq} \displaystyle \mathbb{E}(e^{\sigma|W_{t}|}) \leq 2e^{\sigma^{2}t/2} \end{equation}
The inequality in (\ref{absineq}) can be obtained by direct calculation as follows. \begin{align}\displaystyle \mathbb{E}(e^{\sigma|W_{t}|}) &= \frac{1}{\sqrt{2\pi t}}\left\{\int_{0}^{\infty} e^{(\sigma x - x^{2}/2t)}dx + \int_{-\infty}^{0} e^{-(\sigma x + x^{2}/2t)}dx \right\} \nonumber \\ &= 2e^{\sigma^{2}t/2}\left(\frac{1}{\sqrt{2\pi t}}\int_{0}^{\infty}e^{-(x-t\sigma)^{2}/2t}dx \right) = 2e^{\sigma^{2}t/2}(1 - \Phi(-\sigma\sqrt{t})) \leq 2e^{\sigma^{2}t/2}, \nonumber \end{align} where $\Phi$ denotes the standard normal distribution function.

Recall that we may write $W_{k} = \rho_{k}U_{k} + \sqrt{1-\rho_{k}^{2}}M_{k}^{(1)}$, where ${M_{k}^{(1)}}$ is a Brownian motion process independent of ${U_{k}}$. With this, we have \begin{align} W_{k} + U_{k} &= (1 + \rho_{k})U_{k} + \sqrt{1 - \rho_{k}^{2}}M_{k}^{(1)} \leq 2|U_{k}| + |M_{k}^{(1)}| \nonumber \\ \mathbb{E}(e^{\sigma(W_{k} + U_{k})}) &\leq Ce^{2\sigma^{2}k}e^{\sigma^{2}k/2} = Ce^{5\sigma^{2}k/2} \nonumber \\ \mathbb{E}(G_{1}^{2}.\mathrm{term}4) &\leq C e^{-2c\sigma^{2}T}\sum_{k=1}^{T}\sum_{j=1}^{T}e^{-b\sigma^{2}(k+j)}e^{\frac{\sigma^{2}j}{2}}e^{\frac{5\sigma^{2}k}{2}} \nonumber \\ &= Ce^{-2c\sigma^{2}T}e^{(3.5-b)\sigma^{2}}\left[\frac{1 - e^{(1-b)\sigma^{2}T}}{1-e^{(1-b)\sigma^{2}}} \right] \left[\frac{1 - e^{(2.5-b)\sigma^{2}T}}{1-e^{(2.5-b)\sigma^{2}} } \right] \nonumber \\ & \rightarrow 0\;\mathrm{as}\; T \rightarrow \infty, \;\mathrm{if}\; c > 0\;\mathrm{and}\; b > 2.5. \nonumber \end{align} Here and below $C$ denotes a generic positive finite constant ($0 < C < \infty$). The same type of calculation for $\mathbb{E}(G_{1}^{2}.\mathrm{term}13)$ also.

When $k \geq j$, \begin{align} \mathbb{E}(e^{W_{k} + U_{j}}) &= \mathbb{E}(e^{\sigma [\rho_{k}U_{k} + U_{j} + \sqrt{1 - \rho_{k}^{2}} M_{k}^{(1)} ] }) \leq   \mathbb{E}(e^{\sigma [|U_{k}| + |U_{j}| + |M_{k}^{(1)}| ] }) \nonumber \\ &= \mathbb{E}(e^{\sigma [|U_{k}| + |U_{j}|] })\mathbb{E}(e^{\sigma [|M_{k}^{(1)}| ] })  \leq  \mathbb{E}(e^{\sigma [|U_{k} - U_{j}|] })\mathbb{E}(e^{2\sigma|U_{j}| })\mathbb{E}(e^{\sigma [|M_{k}^{(1)}| ] }) \nonumber \\ &\leq Ce^{\frac{\sigma^{2}}{2}(k-j)}e^{2\sigma^{2}j}e^{\frac{\sigma^{2}}{2}k} \leq Ce^{\frac{3\sigma^{2}}{2}(k+j)}. \nonumber \end{align} Because of the symmetry in the upper bound, the above inequality holds for $k \leq j$ as well. Hence \begin{align}\displaystyle  \mathbb{E}(G_{1}^{2}.\mathrm{term7}) &\leq e^{-2c\sigma^{2}T} \sum_{k=1}^{T}\sum_{j=1}^{T}e^{-b\sigma^{2}(k+j)+\frac{\sigma^{2}}{2}\sigma^{2}(k+j)+\frac{3
\sigma^{2}}{2}\sigma^{2}(k+j)} = e^{-2c\sigma^{2}T} \sum_{k=1}^{T}\sum_{j=1}^{T}e^{(2-b)\sigma^{2}(k+j)} \nonumber \\ &= e^{-2c\sigma^{2}T}e^{2(2-b)\sigma^{2}}\left[\frac{1 - e^{(2-b)\sigma^{2}T}}{1 - e^{(2-b)\sigma^{2}}} \right]^{2} \rightarrow 0\;\mathrm{as}\; T \rightarrow \infty, \;\mathrm{if}\; c > 0\;\mathrm{and}\; b > 2.  \nonumber \end{align} A similar calculation holds for $\mathbb{E}(G_{1}^{2}.\mathrm{term10})$ also.

\begin{align} \mathbb{E}(e^{\sigma(U_{k} + U_{j})}) &= \mathbb{E}(e^{\sigma(U_{k} - U_{j})})\mathbb{E}(e^{2\sigma U_{j}}) = e^{\frac{\sigma^{2}}{2}(k-j)}e^{2\sigma^{2}j} \leq e^{\frac{3\sigma^{2}}{2}(k+j)}, \nonumber \end{align} when $k \geq j$, but the same holds true for $k \leq j$ also because of symmetry. Hence \begin{align}\displaystyle  \mathbb{E}(G_{1}^{2}.\mathrm{term11}) &= e^{-2c\sigma^{2}T} \sum_{k=1}^{T}\sum_{j=1}^{T}e^{-b\sigma^{2}(k+j)}e^{\frac{\sigma^{2}}{2}(k+j)}\mathbb{E}(e^{\sigma(U_{k} + U_{j})}) \leq e^{-2c\sigma^{2}T} \sum_{k=1}^{T}\sum_{j=1}^{T}e^{-b\sigma^{2}(k+j)}e^{2\sigma^{2}(k+j)} \nonumber \\ &= e^{-2c\sigma^{2}T}e^{2(2-b)\sigma^{2}}\left[\frac{1-e^{(2-b)\sigma^{2}T}}{1-e^{(2-b)\sigma^{2}}} \right]^{2} \rightarrow 0\;\mathrm{as}\; T \rightarrow \infty, \;\mathrm{if}\; c > 0\;\mathrm{and}\; b > 2.   \nonumber \end{align} Same calculation works for $G_{1}^{2}.\mathrm{term6}$ also. So far, we have shown convergence to zero of the expectation the terms 4,6,7,8,10,11,12,13,14,15 and 16 of $G_{1}^{2}$. It remains to show the convergence for terms 1,2,3,5 and 9. The calculations for terms 2,3,5 and 9 are similar. So, we will exhibit just the calculation for $G_{1}^{2}.\mathrm{term5}$ and then for $G_{1}^{2}.\mathrm{term1}$.

\begin{align} W_{k} + W_{j} + U_{j} &= [\rho_{k}U_{k} + \sqrt{1 - \rho_{k}^{2}}M_{k}^{(1)}] + [(1+\rho_{j})U_{j} + \sqrt{1 - \rho_{j}^{2}}M_{j}^{(1)}] \nonumber \\ &\leq |U_{k}| + |M_{k}^{(1)}| + 2|U_{j}| + |M_{j}^{(1)}| \nonumber \\ &\leq |U_{k}-U_{j}| + 3|U_{j}| + |M_{k}^{(1)}-M_{j}^{(1)}| + 2|M_{j}^{(1)}|\;\mathrm{when}\; k \geq j. \nonumber \end{align}
\begin{align}\mathbb{E}(e^{\sigma(W_{k} + W_{j} + U_{j})}) &\leq \mathbb{E}(e^{\sigma|U_{k} - U_{j}|})\mathbb{E}(e^{3\sigma|U_{j}|})\mathbb{E}(e^{\sigma|M_{k}^{(1)} - M_{j}^{(1)}|})\mathbb{E}(e^{2\sigma|M_{j}^{(1)}|}) \nonumber \\ &\leq Ce^{\frac{\sigma^{2}}{2}(k-j)}e^{\frac{9\sigma^{2}}{2}j}e^{\frac{\sigma^{2}}{2}(k-j)}e^{2\sigma^{2}j} = Ce^{\sigma^{2}k}e^{\frac{11\sigma^{2}}{2}j} \leq Ce^{\frac{11\sigma^{2}}{2}(k+j)}. \nonumber \end{align} We will get the same upper bound for $k \leq j$ also because of symmetry. Hence, \begin{align}\displaystyle  \mathbb{E}(G_{1}^{2}.\mathrm{term5}) &\leq Ce^{-2c\sigma^{2}T} \sum_{k=1}^{T}\sum_{j=1}^{T}e^{-b\sigma^{2}(k+j)}e^{\frac{\sigma^{2}}{2}k}e^{\frac{11\sigma^{2}}{2}(k+j)} \nonumber \\ &= Ce^{-2c\sigma^{2}T}e^{(11.5-2b)\sigma^{2}}\left[\frac{1-e^{(5.5-b)\sigma^{2}T}}{1-e^{(5.5-b)\sigma^{2}}} \right]\left[\frac{1-e^{(6.5-b)\sigma^{2}T}}{1-e^{(6.5-b)\sigma^{2}}} \right] \nonumber \\ &\rightarrow 0\;\mathrm{as}\; T \rightarrow \infty, \;\mathrm{if}\; c > 0\;\mathrm{and}\; b > 6.5.   \nonumber \end{align} As mentioned above, the calculations for $G_{1}^{2}$ terms 2,3 and 9 are similar to that of term 5. Thus, it remains to show the calculation for $G_{1}^{2}.\mathrm{term1}$ to complete the proof for $G_{1}^{2}$.
\begin{align} W_{k} + U_{k} + W_{j} + U_{j} &= [(1 + \rho_{k})U_{k} + \sqrt{1 - \rho_{k}^{2}}M_{k}^{(1)}] + [(1+\rho_{j})U_{j} + \sqrt{1 - \rho_{j}^{2}}M_{j}^{(1)}] \nonumber \\ &\leq 2|U_{k}| + |M_{k}^{(1)}| + 2|U_{j}| + |M_{j}^{(1)}| \nonumber \\ &\leq 2|U_{k}-U_{j}| + 4|U_{j}| + |M_{k}^{(1)}-M_{j}^{(1)}| + 2|M_{j}^{(1)}|. \nonumber \end{align}
\[ \mathbb{E}(e^{\sigma(W_{k} + U_{k} + W_{j} + U_{j})}) \leq Ce^{2\sigma^{2}(k-j)}e^{8\sigma^{2}j}e^{\frac{\sigma^{2}}{2}(k-j)}e^{2\sigma^{2}j} \leq Ce^{\frac{15\sigma^{2}}{2}(k+j)}. \] The above inequality was obtained assuming $k \geq j$ but the same upper bounds holds for $k \leq j$ as well. Hence
\begin{align}\displaystyle  \mathbb{E}(G_{1}^{2}.\mathrm{term1}) &\leq e^{-2c\sigma^{2}T} \sum_{k=1}^{T}\sum_{j=1}^{T}e^{(7.5-b)\sigma^{2}(k+j)} = e^{-2c\sigma^{2}T}e^{(7.5-b)\sigma^{2}}\left[\frac{1 - e^{(7.5-b)\sigma^{2}T}}{1 - e^{(7.5-b)\sigma^{2}}} \right]^{2} \nonumber \\ &\rightarrow 0\;\mathrm{as}\; T \rightarrow \infty, \;\mathrm{if}\; c > 0\;\mathrm{and}\; b > 7.5.   \nonumber   \end{align} Putting all the calculations together for the 16 terms involved in $G_{1}^{2}$, we see that \[\mathbb{E}(G_{1}^{2}) \rightarrow 0\;\mathrm{as}\; T \rightarrow \infty, \;\mathrm{if}\; c > 0\;\mathrm{and}\; b > 7.5. \]

Now we move onto the calculations for $G_{2}^{2}$, $G_{3}^{2}$ and $G_{4}^{2}$. We will do the calculation for $G_{4}^{2}$ first as it is relatively easy. Also, since $\mathbb{E}(G_{2}^{2}) = \mathbb{E}(G_{3}^{2})$, we will exhibit only the calculation for $\mathbb{E}(G_{2}^{2})$.
\begin{align}\displaystyle  \mathbb{E}(G_{4}^{2}) &\leq e^{-2c\sigma^{2}T} \sum_{k=1}^{T}\sum_{j=1}^{T}e^{a\sigma^{2}(k+j)}\mathbb{E}\left\{(e^{\sigma W_{t} - e^{\sigma^{2}t/2}})^{2}(e^{\sigma U_{t} - e^{\sigma^{2}t/2}})^{2} \right\} \leq Ce^{-2c\sigma^{2}T} \sum_{k=1}^{T}\sum_{j=1}^{T}e^{a\sigma^{2}(k+j)} \nonumber \\ &= C\left[ \frac{e^{a\sigma^{2}}}{1-e^{a\sigma^{2}}} \right]^{2}\left[e^{-c\sigma^{2}T} - e^{(a-c)\sigma^{2}T} \right]^{2} \rightarrow 0\;\mathrm{as}\; T \rightarrow \infty, \;\mathrm{if}\; c > \max\{a,0\}. \nonumber \end{align}

\begin{align} \displaystyle  G_{2}^{2} &= e^{-2c\sigma^{2}T} \sum_{k=1}^{T}\sum_{j=1}^{T}e^{(a-b)\frac{\sigma^{2}}{2}(k+j)}\left\{e^{\sigma(W_{k} + W_{j} + 2U_{t})} - 2e^{\frac{\sigma^{2}}{2}t}e^{\sigma(W_{k} + W_{j} + U_{t})}  \right. \nonumber \\  &+e^{\sigma^{2}t}e^{\sigma(W_{k}+W_{j})} - e^{\frac{\sigma^{2}}{2}j}e^{\sigma(W_{k}+2U_{t})} +2e^{\frac{\sigma^{2}}{2}(j+t)}e^{\sigma(W_{k} + U_{t})}- e^{\sigma^{2}t}e^{\frac{\sigma^{2}}{2}j}e^{\sigma W_{k}} \nonumber \\
& \left.  -e^{\frac{\sigma^{2}}{2}k}e^{\sigma(W_{j} + 2U_{t})} + 2e^{\frac{\sigma^{2}}{2}(k+t)}e^{\sigma(W_{j} + U_{t})} - e^{\sigma^{2}t}e^{\frac{\sigma^{2}}{2}k}e^{\sigma W_{j}} + e^{\frac{\sigma^{2}}{2}(k+j)}(e^{\sigma U_{t}} - e^{\sigma^{2}t/2})^{2} \right\} \nonumber \end{align}

\begin{align} \mathbb{E}(G_{2}^{2}.(\mathrm{last\,term})) &= e^{-2c\sigma^{2}T} \sum_{k=1}^{T}\sum_{j=1}^{T}e^{(a-b)\frac{\sigma^{2}}{2}(k+j)}e^{\frac{\sigma^{2}}{2}(k+j)}\mathbb{E}(e^{\sigma U_{t}} - e^{\sigma^{2}t/2})^{2} \nonumber \\ & \leq C e^{-2c\sigma^{2}T} \sum_{k=1}^{T}\sum_{j=1}^{T}e^{(a-b)\frac{\sigma^{2}}{2}(k+j)}e^{\frac{\sigma^{2}}{2}(k+j)} \nonumber \\  &= C e^{-2c\sigma^{2}T} e^{(a-b+1)\sigma^{2}}\left[\frac{1-e^{(a-b+1)\frac{\sigma^{2}}{2}T}}{1-e^{(a-b+1)\frac{\sigma^{2}}{2}}} \right]^{2} \nonumber \\ &\rightarrow 0, \;\mathrm{as}\; T \rightarrow \infty, \;\mathrm{if}\; c > 0\;\mathrm{and}\; b > a+1. \nonumber \end{align}
Same proof holds for $G_{2}^{2}.\mathrm{term}3$, $G_{2}^{2}.\mathrm{term}6$ and $G_{2}^{2}.\mathrm{term}9$.

For $G_{2}^{2}.\mathrm{term}4$, we start by noting that \[ W_{k} + 2U_{t} \leq |U_{k}| + |M_{k}^{(1)}| + 2|U_{t}| \leq |U_{k} - U_{t}| + 3|U_{t}| + |M_{k}^{(1)}|. \]
\[ \mathbb{E}(e^{\sigma (W_{k} + 2U_{t})}) \leq Ce^{\frac{\sigma^{2}}{2}(k-t)}e^{\frac{9\sigma^{2}}{2}t}e^{\frac{\sigma^{2}}{2}k} \leq Ce^{\sigma^{2}k}, \;\mathrm{when}\; k \geq t,\] \[ W_{k} + 2U_{t} \leq 2|U_{t} - U_{k}| + 3|U_{k}| + |M_{k}^{(1)}|;\;\mathbb{E}(e^{\sigma (W_{k} + 2U_{t})}) \leq Ce^{2\sigma^{2}(t-k)}e^{\frac{9\sigma^{2}}{2}k}e^{\frac{\sigma^{2}}{2}k} \leq Ce^{3\sigma^{2}k}, \;\mathrm{when}\; k \leq t.  \] Combining both cases, $k \geq t$ and $k \leq t$, we get $\mathbb{E}(e^{\sigma (W_{k} + 2U_{t})}) \leq Ce^{3\sigma^{2}k}.$ Hence \begin{align}\mathbb{E}(G_{2}^{2}.\mathrm{term}4) &\leq e^{-2c\sigma^{2}T} \sum_{k=1}^{T}\sum_{j=1}^{T}e^{(a-b)\frac{\sigma^{2}}{2}(k+j)}e^{\frac{\sigma^{2}}{2}j}e^{3\sigma^{2}k} \leq e^{-2c\sigma^{2}T} \sum_{k=1}^{T}\sum_{j=1}^{T}e^{(6 + a-b)\frac{\sigma^{2}}{2}(k+j)} \nonumber \\ &\leq e^{-2c\sigma^{2}T}e^{(6 + a-b)\sigma^{2}}\left[\frac{1 - e^{(6 + a-b)\frac{\sigma^{2}}{2}T}}{1 - e^{(6 + a-b)\frac{\sigma^{2}}{2}}} \right]^{2} \rightarrow 0, \;\mathrm{as}\; T \rightarrow \infty, \;\mathrm{if}\; c > 0\;\mathrm{and}\; b > a+6. \nonumber \end{align} Same type of calculations as above works for $\mathbb{E}(G_{2}^{2}.\mathrm{term}7)$ also.

For $G_{2}^{2}.\mathrm{term}5$, when $k \geq t$, $W_{k} + U_{t} \leq |U_{k} - U_{t}| + 2|U_{t}| + |M_{k}^{(1)}|$, so that \[\mathbb{E}(e^{\sigma(W_{k} + U_{t})}) \leq C\exp \{\frac{\sigma^{2}}{2}[(k-t) + 4t + k ] \} \leq Ce^{\sigma^{2}k} \] and when $k \leq t$, $W_{k} + U_{t} \leq |U_{t} - U_{k}| + 2|U_{k}| + |M_{k}^{(1)}|$, so that \[\mathbb{E}(e^{\sigma(W_{k} + U_{t})}) \leq C\exp \{\frac{\sigma^{2}}{2}[(t-k) + 4k + k ] \} \leq Ce^{\sigma^{2}k}. \] Thus for any $t, k$, $\mathbb{E}(e^{\sigma(W_{k} + U_{t})}) \leq Ce^{\sigma^{2}k}$. Hence \begin{align} \mathbb{E}(G_{2}^{2}.\mathrm{term}5) &\leq e^{-2c\sigma^{2}T} \sum_{k=1}^{T}\sum_{j=1}^{T}e^{(a-b)\frac{\sigma^{2}}{2}(k+j)}e^{\frac{\sigma^{2}}{2}j}e^{2\sigma^{2}k} \nonumber \leq Ce^{-2c\sigma^{2}T}e^{(4+a-b)\sigma^{2}}\left[\frac{1-e^{(4+a-b)\sigma^{2}T/2}}{1-e^{(4+a-b)\sigma^{2}/2}} \right]^{2} \nonumber \\ & \rightarrow 0, \;\mathrm{as}\; T \rightarrow \infty, \;\mathrm{if}\; c > 0\;\mathrm{and}\; b > a+4. \nonumber \end{align} Same type of calculations as above works for $\mathbb{E}(G_{2}^{2}.\mathrm{term}8)$ also.

For $G_{2}^{2}.\mathrm{term}2$, we consider four cases: $\{k \geq t, j \geq t\}$, $\{k \geq t, j \leq t\}$, $\{k \leq t, j \geq t\}$ and $\{k \leq t, j \leq t\}$.

\underline{Case 1a}: $\{k \geq j \geq t\}$. \begin{align} W_{k} + W_{j} + U_{t} & \leq |U_{k}| + |M_{k}^{(1)}| + |U_{j}| + |M_{j}^{(1)}| + |U_{t}| \nonumber \\ & \leq |U_{k} - U_{j}| + 2|U_{j}| + |U_{t}| + |M_{k}^{(1)} - M_{j}^{(1)}| + 2|M_{j}| \nonumber \\ & \leq |U_{k} - U_{j}| + 2|U_{j}-U_{t}| + 3|U_{t}| + |M_{k}^{(1)} - M_{j}^{(1)}| + 2|M_{j}|, \nonumber \end{align} \[\mathbb{E}(e^{\sigma(W_{k} + W_{j} + U_{t})}) \leq C\exp \{\frac{\sigma^{2}}{2}[(k-j) + 4(j-t) + 9t + (k-j)+4j] \} \leq Ce^{\sigma^{2}(k+3j)}. \]

\underline{Case 1b}: $\{j \geq k \geq t\}$. The result for this case can be obtained by interchanging $k$ and $j$ in the result for Case 1a. \[\mathbb{E}(e^{\sigma(W_{k} + W_{j} + U_{t})}) \leq Ce^{\sigma^{2}(3k+j)}. \]

\underline{Case 2}: $\{k \geq t \geq j\}$. \begin{align} W_{k} + W_{j} + U_{t} & \leq |U_{k} - U_{t}| + 2|U_{t}| + |U_{j}| + |M_{k}^{(1)} - M_{t}^{(1)}| + |M_{t}^{(1)}| + |M_{j}^{(1)}| \nonumber \\ & \leq |U_{k} - U_{t}| + 2|U_{t} - U_{j}| + 3|U_{j}| + |M_{k}^{(1)} - M_{t}^{(1)}| + |M_{t}^{(1)}-M_{j}^{(1)}| + 2|M_{j}^{(1)}|, \nonumber \end{align} \[\mathbb{E}(e^{\sigma(W_{k} + W_{j} + U_{t})}) \leq C\exp \{\frac{\sigma^{2}}{2}[(k-t) + 4(t-j) + 9j + (k-t)+(t-j)+4j] \} \leq Ce^{\sigma^{2}(k+4j)}. \]

\underline{Case 3}: $\{j \geq t \geq k\}$. The result for this case can be obtained by interchanging $k$ and $j$ in the result for Case 2. \[\mathbb{E}(e^{\sigma(W_{k} + W_{j} + U_{t})}) \leq Ce^{\sigma^{2}(4k+j)}. \]

\underline{Case 4a}: $\{k \leq j \leq t\}$. \[W_{k} + W_{j} + U_{t} \leq 3|U_{k}| + 2|M_{k}^{(1)}| + 2|U_{j}-U_{k}| + |M_{j}^{(1)}-M_{k}^{(1)}| + |U_{t}-U_{j}|. \] \[\mathbb{E}(e^{\sigma(W_{k} + W_{j} + U_{t})}) \leq C\exp \{\frac{\sigma^{2}}{2}[9k + 4k + 4(j-k) + (j-k) + (t-j)] \} \leq Ce^{\sigma^{2}(4k+2j)}. \]

\underline{Case 4b}: $\{j \leq k \leq t\}$. By interchanging $k$ and $j$ in the result for Case 4a, we get \[\mathbb{E}(e^{\sigma(W_{k} + W_{j} + U_{t})}) \leq Ce^{\sigma^{2}(2k+4j)}. \]

For all cases combined we see that \[\mathbb{E}(e^{\sigma(W_{k} + W_{j} + U_{t})}) \leq Ce^{4\sigma^{2}(k+j)}. \] Hence \begin{align} \mathbb{E}(G_{2}^{2}.\mathrm{term}2) &\leq e^{-2c\sigma^{2}T} \sum_{k=1}^{T}\sum_{j=1}^{T}e^{(a-b)\frac{\sigma^{2}}{2}(k+j)}e^{4\sigma^{2}(k+j)} \nonumber \leq Ce^{-2c\sigma^{2}T}e^{(8+a-b)\sigma^{2}}\left[\frac{1-e^{(8+a-b)\sigma^{2}T/2}}{1-e^{(8+a-b)\sigma^{2}/2}} \right]^{2} \nonumber \\ & \rightarrow 0, \;\mathrm{as}\; T \rightarrow \infty, \;\mathrm{if}\; c > 0\;\mathrm{and}\; b > a+8. \nonumber \end{align}

For $G_{2}^{2}.\mathrm{term}1$ also we again deal with four cases.

\underline{Case 1a}: $\{k \geq j \geq t\}$. \begin{align} W_{k} + W_{j} + 2U_{t} & \leq |U_{k}| + |M_{k}^{(1)}| + |U_{j}| + |M_{j}^{(1)}| + 2|U_{t}| \nonumber \\ & \leq |U_{k} - U_{j}| + 2|U_{j}| + 2|U_{t}| + |M_{k}^{(1)} - M_{j}^{(1)}| + 2|M_{j}| \nonumber \\ & \leq |U_{k} - U_{j}| + 2|U_{j}-U_{t}| + 4|U_{t}| + |M_{k}^{(1)} - M_{j}^{(1)}| + 2|M_{j}|, \nonumber \end{align} \[\mathbb{E}(e^{\sigma(W_{k} + W_{j} + 2U_{t})}) \leq C\exp \{\frac{\sigma^{2}}{2}[(k-j) + 4(j-t) + 16t + (k-j)+4j] \} \leq Ce^{\frac{\sigma^{2}}{2}(2k+6j)}. \]

\underline{Case 1b}: $\{j \geq k \geq t\}$. By interchanging $k$ and $j$ in the result for Case 1a, we get \[\mathbb{E}(e^{\sigma(W_{k} + W_{j} + 2U_{t})}) \leq Ce^{\frac{\sigma^{2}}{2}(6k+2j)}. \]

\underline{Case 2}: $\{k \geq t, j \leq t\}$. \[ W_{k} + W_{j} + 2U_{t} \leq |U_{k} - U_{t}| + 3|U_{t} - U_{j}| + 4|U_{j}| + |M_{k}^{(1)} - M_{t}^{(1)}| + |M_{t}^{(1)}-M_{j}^{(1)}| + 2|M_{j}^{(1)}| \] \[\mathbb{E}(e^{\sigma(W_{k} + W_{j} + 2U_{t})}) \leq Ce^{\frac{\sigma^{2}}{2}(2k+10j)}. \]

\underline{Case 3}: $\{k \leq t, j \geq t\}$. By interchanging $k$ and $j$ in Case 2, we get, \[\mathbb{E}(e^{\sigma(W_{k} + W_{j} + 2U_{t})}) \leq Ce^{\frac{\sigma^{2}}{2}(10k+2j)}. \]

\underline{Case 4a}: $\{k \leq j \leq t\}$. \[W_{k} + W_{j} + 2U_{t} \leq 2|U_{t}-U_{j}| + 3|U_{j}-U_{k}| + 4|U_{k}| + |M_{j}^{(1)} - M_{k}^{(1)}| + 2|M_{k}^{(1)}| \] \[\mathbb{E}(e^{\sigma(W_{k} + W_{j} + 2U_{t})}) \leq Ce^{\frac{\sigma^{2}}{2}(10k+6j)}. \]

\underline{Case 4b}: \[\mathbb{E}(e^{\sigma(W_{k} + W_{j} + 2U_{t})}) \leq Ce^{\frac{\sigma^{2}}{2}(6k+10j)}. \]

Combining all cases, \begin{align} \mathbb{E}(G_{2}^{2}.\mathrm{term}1) &\leq e^{-2c\sigma^{2}T} \sum_{k=1}^{T}\sum_{j=1}^{T}e^{(a-b)\frac{\sigma^{2}}{2}(k+j)}e^{5\sigma^{2}(k+j)} \nonumber \leq Ce^{-2c\sigma^{2}T}e^{(10+a-b)\sigma^{2}}\left[\frac{1-e^{(10+a-b)\sigma^{2}T/2}}{1-e^{(10+a-b)\sigma^{2}/2}} \right]^{2} \nonumber \\ & \rightarrow 0, \;\mathrm{as}\; T \rightarrow \infty, \;\mathrm{if}\; c > 0\;\mathrm{and}\; b > a+10. \nonumber \end{align}

For the remaining terms $G_{1}G_{2}, G_{1}G_{3}, \ldots, G_{3}G_{4}$ we work out the calculations for only the first term in each of them. The conditions on $a$, $b$ and $c$ required for convergence of the first term in each product, will ensure the convergence for the remaining terms in each as well. We start with $G_{1}G_{2}$. \[ \displaystyle  G_{1}G_{2} = e^{-2c\sigma^{2}T} \sum_{k=1}^{T}\sum_{j=1}^{T}e^{-b\sigma^{2}k}e^{(a-b)\frac{\sigma^{2}}{2}j}\{(e^{\sigma W_{k}} - e^{\frac{\sigma^{2}k}{2}})(e^{\sigma U_{k}} - e^{\frac{\sigma^{2}k}{2}})(e^{\sigma W_{j}} - e^{\frac{\sigma^{2}j}{2}})(e^{\sigma U_{t}} - e^{\frac{\sigma^{2}t}{2}}) \} \] \[ \displaystyle  G_{1}G_{2}.\mathrm{term}1 = e^{-2c\sigma^{2}T} \sum_{k=1}^{T}\sum_{j=1}^{T}e^{-b\sigma^{2}k}e^{(a-b)\frac{\sigma^{2}}{2}j}\left[e^{\sigma (W_{k} + U_{k} + W_{j} + U_{t})}\right]. \] \begin{align}W_{k} + U_{k} + W_{j} + U_{t} &= (1+\rho_{k})U_{k} + \sqrt{1-\rho_{k}^{2}}M_{k}^{(1)} + \rho_{j}U_{j} + \sqrt{1-\rho_{j}^{2}}M_{j}^{(1)} + U_{t} \nonumber \\ & \leq 2|U_{k}| + |M_{k}^{(1)}| + |U_{j}| + |M_{j}^{(1)}| + |U_{t}| \nonumber \end{align}

\underline{Case 1a}: $\{k \geq j \geq t\}$.  \[ W_{k} + U_{k} + W_{j} + U_{t} \leq 2|U_{k} - U_{j}| + 3|U_{j} - U_{t}| + 4|U_{t}| + |M_{k}^{(1)} - M_{j}^{(1)}| + 2|M_{j}^{(1)}|. \] \[\mathbb{E}(e^{\sigma (W_{k} + U_{k} + W_{j} + U_{t})}) \leq Ce^{\frac{\sigma^{2}}{2}(5k + 8j)}. \]

\underline{Case 1b}: $\{j \geq k \geq t\}$.  \[ W_{k} + U_{k} + W_{j} + U_{t} \leq 2|U_{j} - U_{k}| + 3|U_{k} - U_{t}| + 4|U_{t}| + |M_{j}^{(1)} - M_{k}^{(1)}| + 2|M_{k}^{(1)}|. \] \[\mathbb{E}(e^{\sigma (W_{k} + U_{k} + W_{j} + U_{t})}) \leq Ce^{\frac{\sigma^{2}}{2}(11k + 2j)}. \]

\underline{Case 2}: $\{k \geq t \geq j\}$.  \[ W_{k} + U_{k} + W_{j} + U_{t} \leq 2|U_{k} - U_{t}| + 3|U_{t} - U_{j}| + 4|U_{j}| + |M_{k}^{(1)} - M_{j}^{(1)}| + 2|M_{j}^{(1)}|. \] \[\mathbb{E}(e^{\sigma (W_{k} + U_{k} + W_{j} + U_{t})}) \leq Ce^{\frac{\sigma^{2}}{2}(5k + 10j)}. \]

\underline{Case 3}: $\{j \geq t \geq k\}$.  \[ W_{k} + U_{k} + W_{j} + U_{t} \leq |U_{j} - U_{t}| + 2|U_{t} - U_{k}| + 4|U_{k}| + |M_{j}^{(1)} - M_{k}^{(1)}| + 2|M_{k}^{(1)}|. \] \[\mathbb{E}(e^{\sigma (W_{k} + U_{k} + W_{j} + U_{t})}) \leq Ce^{\frac{\sigma^{2}}{2}(15k + 2j)}. \]

\underline{Case 4a}: $\{t \geq k \geq j\}$.  \[ W_{k} + U_{k} + W_{j} + U_{t} \leq |U_{t} - U_{k}| + 3|U_{k} - U_{j}| + 4|U_{j}| + |M_{k}^{(1)} - M_{j}^{(1)}| + 2|M_{j}^{(1)}|. \] \[\mathbb{E}(e^{\sigma (W_{k} + U_{k} + W_{j} + U_{t})}) \leq Ce^{\frac{\sigma^{2}}{2}(9k + 10j)}. \]

\underline{Case 4b}: $\{t \geq j \geq k\}$.  \[ W_{k} + U_{k} + W_{j} + U_{t} \leq |U_{t} - U_{j}| + 2|U_{j} - U_{k}| + 4|U_{k}| + |M_{j}^{(1)} - M_{k}^{(1)}| + 2|M_{k}^{(1)}|. \] \[\mathbb{E}(e^{\sigma (W_{k} + U_{k} + W_{j} + U_{t})}) \leq Ce^{\frac{\sigma^{2}}{2}(15k + 4j)}. \]

Combining all cases, we get $\mathbb{E}(e^{\sigma (W_{k} + U_{k} + W_{j} + U_{t})}) \leq Ce^{\frac{\sigma^{2}}{2}(15k + 10j)}. $
\begin{align} \mathbb{E}(G_{1}G_{2}.\mathrm{term}1) & \leq Ce^{-2c\sigma^{2}T} \sum_{k=1}^{T}e^{-b\sigma^{2}k}e^{15k\frac{\sigma^{2}}{2}} \left\{ \sum_{j=1}^{T}e^{(a-b)\frac{\sigma^{2}}{2}j}e^{10j\frac{\sigma^{2}}{2}} \right\} \nonumber \\ & \leq Ce^{-2c\sigma^{2}T}e^{(25+a-3b)\frac{\sigma^{2}}{2}}\left[\frac{1-e^{(10+a-b)\sigma^{2}T/2}}{1-e^{(10+a-b)\sigma^{2}/2}} \right]\left[\frac{1-e^{(15-2b)\sigma^{2}T/2}}{1-e^{(15-2b)\sigma^{2}/2}} \right] \nonumber \\ & \rightarrow 0, \;\mathrm{as}\; T \rightarrow \infty, \;\mathrm{if}\; c > 0\;\mathrm{and}\; b > \max\{a+10,15/2\} = a + 10\; (\mathrm{since\;} a > 0). \nonumber\end{align}

We skip the proof for $G_{1}G_{3}$, since it is very similar to that of $G_{1}G_{2}$ (and the results are the same). We move onto $G_{1}G_{4}$. \[ \displaystyle  G_{1}G_{4} = e^{-2c\sigma^{2}T} \sum_{k=1}^{T}\sum_{j=1}^{T}e^{-b\sigma^{2}k}e^{a\sigma^{2}j}\{(e^{\sigma W_{k}} - e^{\frac{\sigma^{2}k}{2}})(e^{\sigma U_{k}} - e^{\frac{\sigma^{2}k}{2}})(e^{\sigma W_{t}} - e^{\frac{\sigma^{2}t}{2}})(e^{\sigma U_{t}} - e^{\frac{\sigma^{2}t}{2}}) \} \] \[ \displaystyle  G_{1}G_{4}.\mathrm{term}1 = e^{-2c\sigma^{2}T} \sum_{k=1}^{T}\sum_{j=1}^{T}e^{-b\sigma^{2}k}e^{a\sigma^{2}j}\left[e^{\sigma (W_{k} + U_{k} + W_{t} + U_{t})}\right]. \] \[ W_{k} + U_{k} + W_{t} + U_{t} \leq 2|U_{k}| + |M_{k}^{(1)}| + 2|U_{t}| + |M_{t}^{(1)}|. \]  When $k \geq t$, \[W_{k} + U_{k} + W_{t} + U_{t} \leq 2|U_{k} - U_{t}| + 4|U_{t}| + |M_{k}^{(1)} - M_{t}^{(1)}| + 2|M_{t}^{(1)}| \] and when $k \leq t$, \[W_{k} + U_{k} + W_{t} + U_{t} \leq 2|U_{t} - U_{k}| + 4|U_{k}| + |M_{t}^{(1)} - M_{k}^{(1)}| + 2|M_{k}^{(1)}|.  \] In the first case, we will get $\mathbb{E}(e^{\sigma (W_{k} + U_{k} + W_{t} + U_{t})}) \leq Ce^{\frac{5k\sigma^{2}}{2}}$ and in the second case we will get $\mathbb{E}(e^{\sigma (W_{k} + U_{k} + W_{t} + U_{t})}) \leq Ce^{\frac{15k\sigma^{2}}{2}}$ and so combining both cases we see that $\mathbb{E}(e^{\sigma (W_{k} + U_{k} + W_{t} + U_{t})}) \leq Ce^{\frac{15k\sigma^{2}}{2}}$. Hence, \begin{align} \mathbb{E}(G_{1}G_{4}.\mathrm{term}1) & \leq Ce^{-2c\sigma^{2}T}e^{a\sigma^{2}}\left[\frac{1-e^{a\sigma^{2}T}}{1-e^{a\sigma^{2}}} \right]e^{(15-2b)\frac{\sigma^{2}}{2}}\left[\frac{1-e^{(15-2b)\sigma^{2}T/2}}{1-e^{(15-2b)\sigma^{2}/2}} \right] \nonumber \\ & = C e^{(15-2b)\frac{\sigma^{2}}{2}}\left(\frac{e^{a\sigma^{2}}}{1-e^{a\sigma^{2}}} \right)\left[\frac{1-e^{(15-2b)\sigma^{2}T/2}}{1-e^{(15-2b)\sigma^{2}/2}} \right][e^{-2c\sigma^{2}T} - e^{(a-2c)\sigma^{2}T} ] \nonumber \\ &\rightarrow 0, \;\mathrm{as}\; T \rightarrow \infty, \;\mathrm{if}\; b > 15/2\;\mathrm{and}\; c > \max\{0,a/2\} = a/2\; (\mathrm{since\;} a > 0). \nonumber \end{align}

\[ \displaystyle  G_{2}G_{3} = e^{-2c\sigma^{2}T} \sum_{k=1}^{T}\sum_{j=1}^{T}e^{(a-b)\frac{\sigma^{2}}{2}(k+j)}\{(e^{\sigma W_{k}} - e^{\frac{\sigma^{2}k}{2}})(e^{\sigma U_{t}} - e^{\frac{\sigma^{2}t}{2}})(e^{\sigma W_{t}} - e^{\frac{\sigma^{2}t}{2}})(e^{\sigma U_{j}} - e^{\frac{\sigma^{2}j}{2}}) \} \] \[ \displaystyle  G_{2}G_{3}.\mathrm{term}1 = e^{-2c\sigma^{2}T} \sum_{k=1}^{T}\sum_{j=1}^{T}e^{(a-b)\frac{\sigma^{2}}{2}(k+j)}\left[e^{\sigma (W_{k} + U_{j} + W_{t} + U_{t})}\right]. \] \[ W_{k} + U_{j} + W_{t} + U_{t} \leq |U_{k}| + |M_{k}^{(1)}| + |U_{j}| + 2|U_{t}| + |M_{t}^{(1)}| \]

\underline{Case 1a}: $\{k \geq j \geq t\}$.  \[ W_{k} + U_{j} + W_{t} + U_{t} \leq |U_{k} - U_{j}| + 2|U_{j} - U_{t}| + 4|U_{t}| + |M_{k}^{(1)} - M_{t}^{(1)}| + 2|M_{t}^{(1)}|. \] \[\mathbb{E}(e^{\sigma (W_{k} + U_{j} + W_{t} + U_{t})}) \leq Ce^{\frac{\sigma^{2}}{2}(2k + 3j)}. \]

\underline{Case 1b}: $\{j \geq k \geq t\}$.  \[ W_{k} + U_{j} + W_{t} + U_{t} \leq |U_{j} - U_{k}| + 2|U_{k} - U_{t}| + 4|U_{t}| + |M_{k}^{(1)} - M_{t}^{(1)}| + 2|M_{t}^{(1)}|. \] \[\mathbb{E}(e^{\sigma (W_{k} + U_{j} + W_{t} + U_{t})}) \leq Ce^{\frac{\sigma^{2}}{2}(4k + j)}. \]

\underline{Case 2}: $\{k \geq t \geq j\}$.  \[ W_{k} + U_{j} + W_{t} + U_{t} \leq |U_{k} - U_{j}| + 3|U_{t} - U_{j}| + 4|U_{j}| + |M_{k}^{(1)} - M_{t}^{(1)}| + 2|M_{t}^{(1)}|. \] \[\mathbb{E}(e^{\sigma (W_{k} + U_{j} + W_{t} + U_{t})}) \leq Ce^{\frac{\sigma^{2}}{2}(2k + 7j)}. \]

\underline{Case 3}: $\{j \geq t \geq k\}$.  \[ W_{k} + U_{j} + W_{t} + U_{t} \leq |U_{j} - U_{t}| + 3|U_{t} - U_{k}| + 4|U_{k}| + |M_{t}^{(1)} - M_{k}^{(1)}| + 2|M_{k}^{(1)}|. \] \[\mathbb{E}(e^{\sigma (W_{k} + U_{j} + W_{t} + U_{t})}) \leq Ce^{\frac{\sigma^{2}}{2}(10k + j)}. \]

\underline{Case 4a}: $\{k \leq j \leq t\}$.  \[ W_{k} + U_{j} + W_{t} + U_{t} \leq 2|U_{t} - U_{j}| + 3|U_{j} - U_{k}| + 4|U_{k}| + |M_{t}^{(1)} - M_{k}^{(1)}| + 2|M_{k}^{(1)}|. \] \[\mathbb{E}(e^{\sigma (W_{k} + U_{j} + W_{t} + U_{t})}) \leq Ce^{\frac{\sigma^{2}}{2}(10k + 5j)}. \]

\underline{Case 4b}: $\{j \leq k \leq t\}$.  \[ W_{k} + U_{j} + W_{t} + U_{t} \leq 2|U_{t} - U_{k}| + 3|U_{k} - U_{j}| + 4|U_{j}| + |M_{t}^{(1)} - M_{k}^{(1)}| + 2|M_{k}^{(1)}|. \] \[\mathbb{E}(e^{\sigma (W_{k} + U_{j} + W_{t} + U_{t})}) \leq Ce^{\frac{\sigma^{2}}{2}(8k + 7j)}. \]

Overall, the following upper bound is satisfied in all the cases: $\mathbb{E}(e^{\sigma (W_{k} + U_{j} + W_{t} + U_{t})}) \leq Ce^{\frac{\sigma^{2}}{2}[10(k + j)]}.$ Hence we have, \begin{align} \mathbb{E}(G_{2}G_{3}.\mathrm{term}1) & \leq Ce^{-2c\sigma^{2}T}\sum_{k=1}^{T}\sum_{j=1}^{T}e^{(10 + a-b)\frac{\sigma^{2}}{2}(k+j)} \leq Ce^{-2c\sigma^{2}T}e^{(10+a-b)\frac{\sigma^{2}}{2}}\left[\frac{1 - e^{(10+a-b)\frac{\sigma^{2}}{2}T}}{1 - e^{(10+a-b)\frac{\sigma^{2}}{2}}} \right]^{2} \nonumber \\ &\rightarrow 0, \;\mathrm{as}\; T \rightarrow \infty, \;\mathrm{if}\; c > 0\; \mathrm{and}\; b > a + 10. \nonumber   \end{align}

Next we consider $G_{2}G_{4}$. \[ \displaystyle  G_{2}G_{4} = e^{-2c\sigma^{2}T} \sum_{k=1}^{T}\sum_{j=1}^{T}e^{(a-b)\frac{\sigma^{2}}{2}k}e^{a\frac{\sigma^{2}}{2}j}\{(e^{\sigma W_{k}} - e^{\frac{\sigma^{2}k}{2}})(e^{\sigma U_{t}} - e^{\frac{\sigma^{2}t}{2}})(e^{\sigma W_{t}} - e^{\frac{\sigma^{2}t}{2}})(e^{\sigma U_{t}} - e^{\frac{\sigma^{2}t}{2}}) \} \] \[ \displaystyle  G_{2}G_{4}.\mathrm{term}1 = e^{-2c\sigma^{2}T} \sum_{k=1}^{T}\sum_{j=1}^{T}e^{(a-b)\frac{\sigma^{2}}{2}k}e^{a\frac{\sigma^{2}}{2}j}\left[e^{\sigma (W_{k} + 2U_{t} + W_{t})}\right]. \] \[ W_{k} + 2U_{t} + W_{t} \leq |U_{k}| + |M_{k}^{(1)}| + 3|U_{t}| + |M_{t}^{(1)}| \]

$k \geq t$: \[ W_{k} + 2U_{t} + W_{t} \leq |U_{k} - U_{t}| + 4|U_{t}| + |M_{k}^{(1)} - M_{t}^{(1)}| + 2|M_{t}^{(1)}|. \] \[\mathbb{E}(e^{\sigma (W_{k} + 2U_{t} + W_{t})}) \leq C\exp \{\frac{\sigma^{2}}{2}[(k-t) + 16t + (k-t)  + 4t] \} \leq  Ce^{\sigma^{2}k}. \]

$t \geq k$: \[ W_{k} + 2U_{t} + W_{t} \leq 3|U_{t} - U_{k}| + 4|U_{k}| + |M_{t}^{(1)} - M_{k}^{(1)}| + 2|M_{k}^{(1)}|. \] \[\mathbb{E}(e^{\sigma (W_{k} + 2U_{t} + W_{t})}) \leq C\exp \{\frac{\sigma^{2}}{2}[9(t-k) + 16k + (t-k)  + 4k] \} \leq  Ce^{5\sigma^{2}k}. \]

Combining the the two cases $k \geq t$ and $k \leq t$, we get $\mathbb{E}(e^{\sigma (W_{k} + 2U_{t} + W_{t})}) \leq  Ce^{5\sigma^{2}k}$ so that \begin{align} \mathbb{E}(G_{2}G_{4}.\mathrm{term}1) & \leq Ce^{-2c\sigma^{2}T}\sum_{j=1}^{T}e^{a\sigma^{2}j}\sum_{k=1}^{T}e^{(10 + a - b)\frac{\sigma^{2}}{2}k} \nonumber \\ &\leq Ce^{-2c\sigma^{2}T}e^{(10+a-b)\frac{\sigma^{2}}{2}}\left[\frac{1 - e^{(10+a-b)\frac{\sigma^{2}}{2}T}}{1 - e^{(10+a-b)\frac{\sigma^{2}}{2}}} \right]e^{a\sigma^{2}}\left[\frac{1 - e^{a\sigma^{2}T}}{1 - e^{a\sigma^{2}}} \right] \nonumber \\ &\rightarrow 0, \;\mathrm{as}\; T \rightarrow \infty, \;\mathrm{if}\;  b > a + 10\; \mathrm{and}\; c > \max\{a/2,0\} (= a/2,\; \mathrm{since}\; a > 0). \nonumber   \end{align}

Proof and results for $G_{3}G_{4}$ are very similar to that of $G_{2}G_{4}$ and so we skip it.

Combing through all the terms we see that if $c > a > 0$ and $b > a + 10$ each of the expectations converges to zero (and hence $\mathbb{E} (\hat{\gamma}_{t}^{2}) \rightarrow 0$) as $T \rightarrow \infty$, thereby proving the lemma.

\end{document}